\documentclass[preprint,10pt]{elsarticle}

\usepackage{amsmath,amssymb,mathrsfs}
\usepackage{upgreek}
\usepackage[greek,british]{babel}
\usepackage{hyperref}
\usepackage{mathtools}
\usepackage{caption}
\allowdisplaybreaks 
\usepackage{colortbl,enumitem}
\usepackage[ruled]{algorithm2e}
\usepackage{subcaption}
\usepackage{color}
\usepackage[table,dvipsnames]{xcolor} 
\usepackage{threeparttable} 
\newcommand{\headrow}{\rowcolor{black!20}} 
\usepackage{multirow}
\newcommand{\splitcellc}[2][c]{%
  \begin{tabular}[c]{@{}c@{}}\strut#2\strut\end{tabular}}%

\journal{arXiv}

\begin{document}
\begin{frontmatter}
\title{A stable poro-mechanical formulation for Material Point Methods leveraging overlapping meshes and multi-field ghost penalisation}

\author[du]{Giuliano Pretti} 
\author[du]
{Robert E.~Bird} 
\author[du]{Nathan D. Gavin}
\author[du]{William M.~Coombs} 
\author[du]{Charles E.~Augarde\corref{cor1}} 


\cortext[cor1]{Corresponding author:
charles.augarde@durham.ac.uk}
\address[du]{Department of Engineering, Durham University\\
 Science Site, South Road, Durham, DH1 3LE, UK.\\
}

\begin{abstract}
The Material Point Method (MPM) is widely used to analyse coupled (solid-water) problems under large deformations/displacements. However, if not addressed carefully, MPM u-p formulations for poro-mechanics can be affected by two major sources of instability. Firstly, inf-sup condition violation can arise when the spaces for the displacement and pressure fields are not chosen correctly, resulting in an unstable pressure field. Secondly, the intrinsic nature of particle-based discretisation makes the MPM an unfitted mesh-based method, which can affect the system's condition number and solvability, particularly when background mesh elements are poorly populated.
This work proposes a solution to both problems. The inf-sup condition is avoided using two overlapping meshes, a coarser one for the pressure and a finer one for the displacement. This approach does not require stabilisation of the primary equations since it is stable by design and is particularly valuable for low-order shape functions. As for the system's poor condition number, a face ghost penalisation method is added to both the primary equations, which constitutes a novelty in the context of MPM mixed formulations.
This study frequently makes use of the theories of functional analysis or the unfitted Finite Element Method (FEM). Although these theories may not directly apply to the MPM, they provide a robust and logical basis for the research. These rationales are further supported by three numerical examples, which encompass both elastic and elasto-plastic cases and drained and undrained conditions.


\end{abstract}
\end{frontmatter}

\section{Introduction}
Since its initial publications~\cite{sulsky1994particle,sulsky1995application}, the Material Point Method (MPM) has become a widely used method for modelling solid materials undergoing extreme deformations while maintaining a Lagrangian description of the equations.
This versatility has been tested in several engineering applications, including snow avalanches~\cite{gaume2019investigating,gaume2021mechanisms}, ice dynamics~\cite{sulsky2007using,huth2021generalized}, slope stability~\cite{xie2024stabilised,gavin2024implementation}, soft robots~\cite{davy2023framework}, and biomedical applications~\cite{li2021novel}, to name a few. 
However, the MPM has also been appreciated beyond its engineering purposes, especially in computer graphics simulations~\cite{stomakhin2013material,chen2021momentum}. 

Among the different applications, fields such as geotechnics and biomechanics deal with porous solid materials, whose mechanical behaviour is strongly influenced by the presence of an interstitial fluid. 
These materials are the subject of the study of poromechanics~\cite{coussy2004poromechanics}.
MPM-based analyses investigating the behaviour of these materials have dramatically increased in the last years, as the recent review paper by Zheng \emph{et al.}~\cite{zheng2023material} demonstrates.
The choices available to run these studies are numerous, and they can be grouped into two categories (see, in this regards, Soga \emph{et al.}~\cite{soga2016trends}). 
The use of different primary equations distinguishes between $\boldsymbol{u}-p^{(f)}$, $\boldsymbol{v}-\boldsymbol{w}$, and $\boldsymbol{u}-p^{(f)}-\boldsymbol{w}$ formulations, where these formulations are labelled after the name of the primary unknowns ($\boldsymbol{u}$ or $\boldsymbol{v}$ is the displacement or velocity of the solid body, $p^{(f)}$ is the fluid pressure, and $\boldsymbol{w}$ the fluid velocity). 
Since MPM formulations are based on a point-based discretisation (these are indeed the Material Points, MPs), it is possible to choose how many sets of MPs can be employed in a poro-mechanical simulation. 
In one case, one set of MPs can be considered, simultaneously keeping track of both phases of the mixed body. 
In the other, the solid and the fluid phases are separately represented by two sets of MPs. 
The best application of each combination of primary variables and the number of MP sets firmly depends upon the considered problem, and assessing these is beyond the scope of this work, which adopts a one-set $\boldsymbol{u}-p^{(f)}$ formulation. 
The reader interested in a detailed explanation and classification of the above options can refer to these recent theses and references therein~\cite{nost2019iteratively,zheng2022stabilised,pretti2024continuum}.

One-set $\boldsymbol{u}-p^{(f)}$ formulation of the MPM can present two sources of instabilities following the discretisation process, which can enormously incapacitate the fundamental scope of the analyses. 
The first, exhibited particularly in undrained or nearly-undrained conditions, is due to the combination of spaces for the the displacement and the pressure fields, and is the well-known \emph{inf-sup} or \emph{LBB} (Ladyzhenskaya-Babu\v{s}ka-Brezzi) condition. 
After briefly introducing mixture theory in Section~\ref{sec:mixed_theory}, Section~\ref{sec:continuum_formulation} is devoted to illustrate the necessary requirements to avoid this instability, and in Section~\ref{subsec:Q1-iso-Q2-Q1_discretisation} a straightforward remedy is proposed, which is particularly attractive for low-order interpolating functions. 
A MP-based discretisation (Section~\ref{subsec:MPM_algorithm}) inevitably leads to  small overlaps between the stencils of the shape functions and the physical domain of the MPs (and their role as quadrature points) and this is the source of the second instability (Section~\ref{subsec:small-cut_issue}). 
In the literature, this problem often goes under the name of the \emph{small-cut} issue. 
If the overlap between the shape function's stencil and the physical MP-based domain is practically infinitesimal, the system matrix can become ill-conditioned and the requirements discussed in Section~\ref{sec:continuum_formulation} are not met.
A possible remedy coming from the unfitted FEM literature is therefore evaluated in Section~\ref{subsec:face_ghost_stabilisations}.
It must be noted that these two sources of instabilities are not restricted to the MPM and to a $\boldsymbol{u}-p^{(f)}$ formulation for poro-mechanics: similar findings have been pointed out, for instance, by Burman and Hansbo~\cite{burman2014fictitious} in the context of unfitted FEM for Stokes problems and the applicability of the methodologies discussed in this manuscript go well beyond the MPM.  
 
The lack of a formal analysis (see, for instance, Nguyen \emph{et al.}~\cite{nguyen2023material}) of the MPM inhibits some results, which are pivotal and well-acknowledged in the FEM theory. 
However, the similarities between the two methods can be quite substantial in some aspects. 
This point has been raised by different authors (see, for instance,~\cite{guilkey2003implicit,coombs2020lagrangian}), especially in light of seeing the MPM as a FEM where quadrature points are allowed to move independently from the mesh across different time-steps. 
Exploiting these similarities in some circumstances permits us to refer to functional analysis and the (unfitted) FEM theory and inherit results which, even if not formally applicable, are practically helpful in the MPM context. 
This idea of falling back to the FEM theory for MPM analyses, even if not rigorous, is motivated by rationales provided throughout the paper and supported by numerical examples in Section~\ref{sec:numerical_examples}.

\section{Fundamentals of mixture theory}
\label{sec:mixed_theory}
To formalise the mechanics of a two-phase (solid/fluid) material modelled with one-set of MPs, this section briefly details the underlying continuum formulation for mixture theory.
\subsection{Kinematics of the phases}
The kinematics assumptions considered in this work are as follows:
\begin{enumerate}
\item \label{hyp:two_continua}
the fully saturated porous material is treated as two (solid and fluid) juxtaposed continua;
\item \label{hyp:one_fluid}
for simplicity, only one fluid material constitutes this phase;
\item \label{hyp:no_termal_effects}
thermal effects and viscosity are neglected; and
\item \label{hyp:finite_strain}
the porous material undergoes finite strains and rotations in the elasto-plastic regime.
\end{enumerate}

Having introduced the founding assumptions, let the mixed body $\mathcal{B}$ occupy a volume $\Omega$ in the original configuration (at time $t=0$) and a volume $\omega$ in the current configuration (at the generic time $t$). 
Owing to Assumption~\ref{hyp:two_continua}, each infinitesimal volume belonging to a mixed body contains a solid and a unique (as implied by Assumption~\ref{hyp:one_fluid}) fluid phase. 
These phases are mapped through the different configurations via the invertible mappings $\boldsymbol{\varphi}^{(ph)} \left( \boldsymbol{X}^{(ph)}, t \right),$ with $ph = s,f$ (standing for \emph{s}olid and \emph{f}luid), and $\boldsymbol{X}^{(ph)}$ indicating the initial position of each phase.
These continua can evolve differently in different configurations over time, which raises the question of which phase to trace. 
In this sense, the mapping of the solid configuration is \emph{de facto} privileged over the mapping of the fluid phase, which is mostly undirectly traced by considering the relative velocity between the phases.  
Given this consideration, the solid material particle of initial mixed volume $d \Omega$ occupies a volume $d \omega = J \, d \Omega$ in the current configuration, where the Jacobian $J$ is given by the determinant of the solid phase deformation gradient\footnote{
For the sake of completeness, a second Jacobian relative to the fluid phase could be introduced, this being given by $J^{(f)} = \det \boldsymbol{F}^{(f)}$, with $F^{(f)}_{iI} = \partial \varphi^{(f)}_i / \partial X^{(f)}_I$. 
However, since the mixture theory assumes that the current position and mixed volume are shared by the two phase particles, i.e., $\boldsymbol{x} = \varphi^{(s)} ( \boldsymbol{X}^{(s)}, t ) = \varphi^{(f)} ( \boldsymbol{X}^{(f)}, t )$ and $J = J^{(f)}$, the fluid deformation gradient becomes mostly redundant.
}, this being
\begin{equation}
\boldsymbol{F} \cdot d \boldsymbol{X}^{(s)} = d \boldsymbol{x}, \quad \text{with} \ F_{iI} = \frac{\partial \varphi^{(s)}_{i}}{\partial X^{(s)}_{I}}.
\end{equation}
In the above equation, $\boldsymbol{x}$ denotes the material particle position in the current configuration obtained using the solid mapping (i.e., $\boldsymbol{x} = \boldsymbol{\varphi}^{(s)} ( \boldsymbol{X}^{(s)}, t )$). The indices $i = 1, \dots, 3$ and $I = 1, \dots, 3$ indicate the components of the matching Cartesian basis vectors $\boldsymbol{e}_{i} = \boldsymbol{E}_I$ in the current ($\boldsymbol{e}_{i}$) and reference ($\boldsymbol{E}_{I})$ configurations.

As proposed by Kr{\"o}ner \cite{kroner1959allgemeine}, Lee \cite{lee1969elastic} and Mandel~\cite{mandel1971plasticite} and considered in Assumption~\ref{hyp:finite_strain}, the deformation gradient is multiplicatively decoupled into an elastic and plastic part as follows
\begin{equation}
\boldsymbol{F} = \boldsymbol{F}^e \cdot \boldsymbol{F}^p,
\end{equation}
where , due to Assumption~\ref{hyp:no_termal_effects}, no other effects contribute to the deformation gradient.
Other measures of strain used in this work are based on the deformation gradient (or its elastic/plastic part), such as the \emph{left Cauchy-Green strain} $\boldsymbol{b} \coloneqq \boldsymbol{F} \cdot \boldsymbol{F}^T$  and the logarithmic strain $\boldsymbol{\epsilon} \coloneqq \frac{1}{2} \ln \boldsymbol{b}$.

To describe the volume fractions of each phase, the Eulerian porosity $n$ is defined as the ratio between the current fluid volume of the material particle $d \omega^{(f)}$ and its current mixed volume $d \omega$, i.e., $n \coloneqq \frac{d \omega^{(f)}}{d \omega}$. 
\emph{E contrario}, $\left( 1 - n \right)$ gives the volume fraction of the material particle occupied by the solid phase, i.e., $\left( 1 - n \right) = \frac{d \omega^{(s)}}{d \omega}$. 
The initial value of the Eulerian porosity is denoted by $n_0$.
As can be seen, the Eulerian porosity satisfies the property 
\begin{equation}
\label{ineq:bounded_porosity}
0 < n < 1, 
\end{equation} 
where the extremes are excluded to always consider coexisting phases in the material.

The \emph{displacements} are defined as the difference between the current and the original positions. Given the choice of matching Cartesian basis vectors, displacements can be defined as 
\begin{equation}
\boldsymbol{u}^{(ph)} \coloneqq \boldsymbol{\varphi}^{(ph)} \left( \boldsymbol{X}^{(ph)}, t \right) - \boldsymbol{X}^{(ph)}.
\end{equation}
The total time derivative of the current position gives the \emph{material velocities}
\begin{equation}
\boldsymbol{v}^{(ph)} \coloneqq \frac{d}{dt}\biggm|_{(ph)}  \boldsymbol{\varphi} \left( \boldsymbol{X}^{(ph)}, t \right),\end{equation}
where $\frac{d}{dt}\bigm|_{(ph)} \left( \bullet \right)$ indicates the time derivative of $\left(  \bullet \right)$ following the $ph-$phase.

\subsection{Mass conservation}
On top of the principle of mass conservation, these further assumptions are introduced:
\begin{enumerate}[resume]
    \item \label{hyp:no_mass_exchange}
    the constituents do not exchange mass;
    \item \label{hyp:incompressiblity}
    both constituents are incompressible; and
    \item \label{hyp:homogeneous_porosity}
    the porosity network is homogeneous and connected across the material.
\end{enumerate}
Assumption~\ref{hyp:no_mass_exchange} permits writing the conservation of the mixture separately for each phase, i.e.,
\begin{align}
\label{eq:solid_mass_conservation}
\frac{d}{dt}\biggm|_{(s)} \int_{\omega^{(s)}} \rho^{(s)} \left( 1 - n \right) dv & = 0;
\\
\label{eq:fluid_mass_conservation}
\frac{d}{dt}\biggm|_{(f)} \int_{\omega^{(f)}} \rho^{(f)}  n \, dv & = 0,
\end{align}
where $\rho^{(ph)}$ indicates the \emph{mesoscopic density} of the $ph-$phase, which, according to Assumption~\ref{hyp:incompressiblity}, is constant, i.e., $\rho^{(ph)} = \rho_0^{(ph)}$.
The exclusion of a double porosity network by Assumption~\ref{hyp:homogeneous_porosity} allows the fluid mass balance to be written as in Eq.~\eqref{eq:fluid_mass_conservation}, i.e., considering only one type of porosity, approximatively homogeneous in size and shape in the medium under investigation.
Owing to the arbitrary nature of the integration volume and the incompressibility of the solid phase (part of Assumption~\ref{hyp:incompressiblity}), Eq.~\eqref{eq:solid_mass_conservation} can also be expressed as follows
\begin{equation}
\label{eq:Eulerian_porosity}
n = 1 - \frac{1}{J} \left( 1 - n_0 \right).
\end{equation}
Two other equations are employed to express Eq.~\eqref{eq:fluid_mass_conservation} in a convenient way. The former is given by the relationship between material time derivatives relative to the two different phases (see Thorpe~\cite{thorpe1962momentum}), i.e.,
\begin{equation}
\label{eq:time_der_diff_phases}
\frac{d}{dt}\biggm|_{(f)} \left( \bullet \right) 
= 
\frac{d}{dt}\biggm|_{(s)} \left( \bullet \right)
+
\text{\textbf{grad}} \left( \bullet \right) \cdot \left( \boldsymbol{v}^{(f)} -  \boldsymbol{v}^{(s)} \right),
\end{equation}
where $\text{\textbf{grad}}$ indicates the gradient with respect to the current position $\boldsymbol{x}$.
The second equation stems from the property relative to the derivation of determinants of second-order tensors, and it is given by
\begin{equation}
\label{eq:time_der_Jacobian}
\frac{d}{dt}\biggm|_{(ph)} J = J \, \text{\textbf{div}} \cdot \boldsymbol{v}^{(ph)},
\end{equation}
where $\text{\textbf{div}}$ is the divergence with respect to the current position $\boldsymbol{x}$.
Invoking the arbitrary nature of the integration domain again, and using Eqs.~\eqref{eq:Eulerian_porosity}-\eqref{eq:time_der_Jacobian}, the fluid mass conservation given by Eq.~\eqref{eq:fluid_mass_conservation} can be re-written as
\begin{equation}  
\label{eq:fluid_mass_conservation_strong}
\rho^{(f)}_0 \frac{d J}{dt}\Bigm|_{(s)} + J \, \text{\textbf{div}} \cdot \boldsymbol{q}^{(f)} = 
\rho^{(f)}_0 \dot{J} + J \, \text{\textbf{div}} \cdot \boldsymbol{q}^{(f)} = 0,
\end{equation}
where the shorthand notation $\dot{\left( \bullet \right)} = \frac{d}{dt}\bigm|_{(s)} \left( \bullet \right)$ and $\boldsymbol{q}^{(f)}$ is the relative flux  defined as
\begin{equation}
\boldsymbol{q}^{(f)} \coloneqq \rho^{(f)}_0 \, n \left( \boldsymbol{v}^{(f)} - \boldsymbol{v}^{(s)} \right).
\end{equation}

\subsection{Balance of rate of linear momentum}
Let us further assume that
\begin{enumerate}[resume]
\item \label{hyp:no_inertia}
inertia effects are negligible.
\end{enumerate}
Under this Assumption~\ref{hyp:no_inertia}, the balance of the rate of the linear momentum for the mixed body can be expressed as follows
\begin{equation}
\label{eq:mixed_eq_equation_strong}
\text{\textbf{div}} \cdot \boldsymbol{\sigma}
+
\rho \ \boldsymbol{b} = \boldsymbol{0},
\end{equation}
where $\boldsymbol{\sigma}$ is the  \emph{total Cauchy stress}, $\rho \coloneqq \rho^{(s)}_0 \left( 1 - n \right) +  \rho^{(f)}_0 \, n$ is the \emph{current density of the porous material} and $\boldsymbol{b}$ are the body forces per unit mixed weight.
The stress introduced in the above equation can be further characterised. 
As demonstrated by Borja and Alarc{\'o}n~\cite{borja1995mathematical}, considering a perfect fluid (included in Assumption~\ref{hyp:no_termal_effects}) and an incompressible solid phase (Assumption~\ref{hyp:incompressiblity}) allows the \emph{Terzaghi effective stress decomposition}\footnote{Eq.~\eqref{eq:Terzaghi_decomposition} assumes the continuum mechanics sign convention (positive pressure indicates tensile behaviour) as opposed to the geotechnical one (positive pressure designates compressive behaviour).} to hold even in the finite strain context, this being
\begin{equation}
\label{eq:Terzaghi_decomposition}
\boldsymbol{\sigma} = \boldsymbol{\sigma}' - p^{(f)} \boldsymbol{1}.
\end{equation}
In the above equation and throughout this manuscript, the dash $\left( \bullet \right)'$ denotes the \emph{effective} quantities of stress, i.e., those obeying to a constitutive relationship with the the solid strains $\boldsymbol{\sigma}' (\boldsymbol{u}^{(s)} )$, while $p^{(f)}$ is the \emph{Cauchy fluid pressure}.

\subsection{Constitutive relationships}
To introduce the necessary constitutive relationships, it is assumed that
\begin{enumerate}[resume]
\item \label{hyp:isotropic_medium}
the material is isotropic; and
\item \label{hyp:low_reynolds_no}
the fluid flow exhibits Low Reynolds number;
\end{enumerate}
Assumption~\ref{hyp:isotropic_medium} has two consequences. 
On the one hand, it results in \emph{isotropic} elasto-plastic behaviour relating the effective stress and the solid strains.
In particular, this work adopts the \emph{improved Hencky material} suggested in Pretti \emph{et al.}~\cite{pretti2024preserving} to avoid negative values of the Eulerian porosity. 
The effective free energy function $\Psi'$ of this material is given by
\begin{equation}
\label{eq:strain_energy}
\Psi' \left( \boldsymbol{\epsilon}, \boldsymbol{\alpha} \right) 
= 
\frac{K}{2 \, n} \left( \epsilon_v^e \right)^2 + \frac{3}{2} G \left( \epsilon_q^e \right)^2 
+ 
\tilde{\Psi}' \left( \boldsymbol{\alpha} \right),
\end{equation}
where $K>0$ and $G>0$ are the bulk parameter and the shear modulus, while $\tilde{\Psi}'$ defines a part of the free energy function responsible for the kinematic hardening, based on a set of internal variables $\boldsymbol{\alpha}$. 
A few invariants of the logarithmic strains have also been employed in the above equation, which are defined as 
\begin{equation}
\epsilon_v \coloneqq \boldsymbol{\epsilon} \boldsymbol{:} \boldsymbol{1}; 
\qquad
\boldsymbol{e} \coloneqq \boldsymbol{\epsilon} - \frac{\epsilon_v}{3} \boldsymbol{1};
\qquad
\epsilon_q \coloneqq \sqrt{\frac{2}{3}  \, \boldsymbol{e} \boldsymbol{:} \boldsymbol{e}},
\end{equation}
with $\boldsymbol{:}$ being the double contraction operator between tensors. If isotropic permeability is assumed (second consequence of Assumption~\ref{hyp:homogeneous_porosity}) together with Assumption~\ref{hyp:low_reynolds_no}, the \emph{Darcy-Weisbach} law relates the relative flux to the fluid pressure gradient and body forces as follows
\begin{equation}
\label{eq:darcy}
\boldsymbol{q}^{(f)} = - \frac{\upkappa}{g} \left( \text{\textbf{grad}} \, p^{(f)} - \rho^{(f)} \boldsymbol{b} \right),
\end{equation}
with $\upkappa$ being the \emph{hydraulic conductivity} and $g$ the magnitude of the gravitational force. 
The \emph{Kozeny-Carman} formula is assumed to account for the dependency of the hydraulic conductivity from the porosity, i.e.,
\begin{equation}
\upkappa = c_1 \frac{n^3}{\left( 1 - n \right)^2},
\end{equation}
with $c_1$ being a constant parameter.
As embedded in Assumption~\ref{hyp:incompressiblity}, the fluid pressure is unrelated to any constitutive relationship and, as pointed out by Miehe \emph{et al.}~\cite{miehe2015minimization}, acts as a Lagrange multiplier by shifting the dependency of stresses on strain (i.e., the idea of effective stress) depending on the water mass conservation constraint (i.e., drained/undrained conditions).  
This observation has direct consequences on the mapping processes detailed in Section~\ref{subsec:mapping_processes}, since the fluid pressure does not relative to any form of energy.

\section{Continuum weak form and arising requirements}
\label{sec:continuum_formulation}
This section introduces a weak statement of the equations for a $\boldsymbol{u}-p^{(f)}$ formulation, which results in a saddle point formulation (Section~\ref{subsec:weak_form}). 
As is well-acknowledged, this kind of problem can be solved only under specific conditions.  

\subsection{From strong to weak form of equations}
\label{subsec:weak_form}
The primary equations employed in a $\boldsymbol{u}-p^{(f)}$ formulation are given by the mixture equilibrium Eq.~\eqref{eq:mixed_eq_equation_strong} and the fluid mass conservation Eq.~\eqref{eq:fluid_mass_conservation_strong}, these being
\begin{align}
\label{eq:mixed_eq_equation_strong_1}
\text{\textbf{div}} \cdot \left( \boldsymbol{\sigma}'
- p^{(f)} \, \boldsymbol{1}
\right)
+
\rho \ \boldsymbol{b} & = \boldsymbol{0}
 \quad \text{in}  \ \ \omega;
\\
\label{eq:fluid_mass_conservation_strong_1}
\rho_0^{(f)} \dot{J} - J \,\text{\textbf{div}} \cdot \left( \frac{\upkappa}{g} \left( \text{\textbf{grad}} \, p^{(f)} - \rho^{(f)} \boldsymbol{b} \right) \right) & = 0
\quad \text{in} \ \ \omega^{(f)},
\end{align}
where the Terzaghi effective stress decomposition Eq.~\eqref{eq:Terzaghi_decomposition} and the Darcy-Weisbach law Eq.~\eqref{eq:darcy} have been substituted.
The current boundary $\gamma = \partial \omega$ of the considered mixed body is partitioned as follows
\begin{gather}
\label{eq:boundary_union}
\gamma 
= 
\gamma^{\bar{u}} \cup \gamma^{\bar{t}}
= 
\gamma^{\bar{p}} \cup \gamma^{\bar{q}};
\\
\label{eq:boundary_intersection}
\gamma^{\bar{u}} \cap \gamma^{\bar{t}} 
=
\emptyset 
= 
\gamma^{\bar{p}} \cap \gamma^{\bar{q}},
\end{gather}
where $\gamma^{(\bullet)}$ is a particular portion of the boundary. 
Prescribed boundary conditions (BCs) are given on these parts and are as follows 
\begin{alignat}{2}
\label{eq:applied_displacements}
\boldsymbol{u}  & = \bar{\boldsymbol{u}} 
\qquad & \hbox{on} \ \gamma^{\bar{u}}; 
\\
\label{eq:applied_tractions}
\boldsymbol{\sigma} \, \boldsymbol{n}  & = \bar{\boldsymbol{t}} 
\qquad & \hbox{on} \ \gamma^{\bar{t}}; 
\\
\label{eq:applied_pressure}
p^{(f)} & = \bar{p} 
\qquad & \hbox{on} \ \gamma^{\bar{p}}; 
\\
\label{eq:applied_fluxes}
\boldsymbol{q} \cdot \, \boldsymbol{n}  & = \bar{q}
\qquad & \hbox{on} \ \gamma^{\bar{q}}.
\end{alignat}
The above equations constitute the Dirichlet, Eqs.~\eqref{eq:applied_displacements} and~\eqref{eq:applied_pressure}, and the Neumann, Eqs.~\eqref{eq:applied_tractions} and~\eqref{eq:applied_fluxes}, BCs for the mixture equilibrium~\eqref{eq:mixed_eq_equation_strong_1} and fluid mass conservation~\eqref{eq:fluid_mass_conservation_strong_1}, respectively.

The weak form of the above problem is obtained by introducing test functions belonging to the required function spaces and integrating over the relative domains\footnote{For the fluid mass conservation, it must be noted that, technically, the correct integration domain is the fluid. However, this domain and that of the mixture are related by the Eulerian porosity, i.e., $d \omega^{(f)} = n \, d \omega$. If this (i.e., the porosity) is simplified for all the terms in the fluid mass conservation, the integration over the mixture volume is \emph{de facto} achieved.}.  
By applying the divergence theorem and using Neumann BCs~\eqref{eq:applied_tractions} and~\eqref{eq:applied_fluxes}, the weak form can be stated as follows: find $\boldsymbol{u} \in \mathscr{V}_{\bar{\boldsymbol{u}}}$ and $p^{(f)} \in {\mathscr{W}_{\bar{p}}}$ such that, for $t \in [ 0, \tilde{t} ] $,
\begin{multline}
\label{eq:weak_mixture_equilibrium}
\delta \Pi^{eq} \left( \left( \boldsymbol{u} ;  p^{(f)}  \right) ;  \delta \boldsymbol{w} \right)  
\coloneqq  
\int_{\omega} \text{\textbf{grad}} \left( \delta \boldsymbol{w} \right)  \boldsymbol{:} \left( \boldsymbol{\sigma}'
- p^{(f)} \, \boldsymbol{1} \right) \ d v 
-
\int_{\omega} \rho \, \delta \boldsymbol{w} \cdot  \boldsymbol{b} \ dv 
\\
-
\int_{\gamma^{\bar{t}}} \delta \boldsymbol{w} \cdot \bar{\boldsymbol{t}} \ da
 = 0, 
\qquad \forall \delta\boldsymbol{w} \in \mathscr{V}_{\boldsymbol{0}};
\end{multline}
\begin{multline}
\label{eq:weak_fluid_mass_conservation}
\delta \Pi^{cons} \left( \left( \boldsymbol{u} ;  p^{(f)}  \right) ;  \delta \eta  \right)  \coloneqq 
\int_{\omega} \rho_0^{(f)} \, \delta \eta \frac{\dot{J}}{J} \ dv 
+ 
\int_{\omega} \frac{\upkappa}{g} \,\text{\textbf{grad}} \left( \delta \eta \right) \cdot  \left( \text{\textbf{grad}} \, p^{(f)} - \rho^{(f)} \boldsymbol{b} \right) dv
\\
-
\int_{\gamma^{\bar{q}}} \delta \eta \, \bar{q} \ da  = 0,
\qquad \forall \delta \eta \in \mathscr{W}_{0},
\end{multline}
giving the initial condition $\boldsymbol{u} \left( t = 0 \right) = \boldsymbol{u}_0$. 
The potentials $\Pi^{eq}$ and $\Pi^{cons}$ and their first variation $\delta \left( \bullet \right) \left(  \left( \dots \right) ; \delta \left( \tilde{\bullet} \right) \right) $ with respect to their argument $\left( \tilde{\bullet} \right)$ have been above employed, while the spaces of trial functions are defined as follows
\begin{align}
\mathscr{V}_{\bar{\boldsymbol{u}}} 
& = 
\left\{
\boldsymbol{u} \in \left[ H^{1} \left( \omega \right) \right]^{n^{dim}}
\bigm| 
\boldsymbol{u} = \bar{\boldsymbol{u}} 
\ \hbox{on} \ \gamma^{\bar{u}}
\right\};
\\
\label{eq:space_pressures}
\mathscr{W}_{\bar{p}} 
& =
\left\{
p^{(f)} \in \left[ H^{1} \left( \omega \right) \right]
\bigm| 
p^{(f)}  = \bar{p} 
\ \hbox{on} \ \gamma^{\bar{p}}
\right\},
\end{align}
where $H^{1} \left( \omega \right)$ denotes the Sobolev space of degree one on $\omega$. The spaces of the test functions $\delta \boldsymbol{w}$ and $\delta \eta$ are denoted by $\mathscr{V}_{\boldsymbol{0}}$ and $\mathscr{W}_{0}$ and follow similar definitions as the above, with the exception of being zero on the boundary.
It is well-known (see, for instance, Dortdivanlioglu \emph{et al.}~\cite{dortdivanlioglu2018mixed}) that the above weak problem Eqs.~\eqref{eq:weak_mixture_equilibrium}-\eqref{eq:weak_fluid_mass_conservation} can be given by the stationarity of the following saddle point problem
\begin{equation}
\inf_{\boldsymbol{w} \in \mathscr{V}_{\boldsymbol{0}} } \sup_{q \in \mathscr{W}_{0} } \Pi \left( \boldsymbol{w} ; \eta \right),
\end{equation}
with $\Pi \left( \boldsymbol{w} ; \eta \right) 
= \Pi^{eq} \left( \boldsymbol{w} ; \eta \right) 
- 
\Pi^{cons} \left( \boldsymbol{w} ; \eta \right) $.

\subsection{Linearisation}
\label{subsec:linearisation}
Eqs.~\eqref{eq:weak_mixture_equilibrium} and~\eqref{eq:weak_fluid_mass_conservation} represent a non-linear system of equations in $\boldsymbol{u}$. 
This is primarily due to having set the problem in the context of finite strain mechanics. On top of this, other non-linearities, such as elasto-plastic behaviour or constitutive equation~\eqref{eq:strain_energy}, can be included too. To solve these equations, the Newton-Raphson (NR) method stipulates that linearisation of  Eqs.~\eqref{eq:weak_mixture_equilibrium} and~\eqref{eq:weak_fluid_mass_conservation} at a solution $\bigl( \check{\boldsymbol{u}}; \check{p}^{(f)} \bigr)$ is required, this being 
\begin{align}
\label{eq:linearised_functional}
\nonumber
0 = 
\delta \Pi \Bigl( \bigl( \check{\boldsymbol{u}}; \check{p}^{(f)} \bigr) ; \left( \delta \boldsymbol{w}, \delta \boldsymbol{\eta} \right) \Bigr)
& \approx 
\delta \Pi \left( \bigl( \boldsymbol{u}; p^{(f)} \bigr) ; \left( \delta \boldsymbol{w}, \delta \boldsymbol{\eta} \right) \right)
\\
\nonumber
& \qquad
+
\delta 
\Biggl(  
\biggl( \delta \Pi 
\Bigl( 
\bigl( \boldsymbol{u}; p^{(f)} \bigr) 
; 
\left( \delta \boldsymbol{w}, \delta \boldsymbol{\eta} \right) 
\Bigr) 
\biggr) ; 
\bigl( \delta \boldsymbol{u}, \delta \boldsymbol{p}^{(f)} \bigr) 
\Biggr)
\\
\nonumber
& = 
\delta \Pi^{eq} 
\Bigl( 
\bigl( \boldsymbol{u}; p^{(f)} \bigr) 
; 
\left( \delta \boldsymbol{w} \right) 
\Bigr)
-
\delta \Pi^{cons} 
\Bigl( 
\bigl( \boldsymbol{u}; p^{(f)} \bigr) 
; 
\bigl( \delta \eta \bigr) 
\Bigr)
\\ 
\nonumber
& \qquad 
+
\delta 
\Biggl(  
\biggl( 
\delta \Pi^{eq} 
\Bigl( 
\bigl( \boldsymbol{u}; p^{(f)} \bigr) 
;  \delta \boldsymbol{w} 
\Bigr) 
\biggr) 
; 
\bigl( \delta \boldsymbol{u}, \delta \boldsymbol{p}^{(f)} \bigr) 
\Biggr)
\\ 
& \qquad \qquad
-
\delta 
\Biggl(
\biggl( 
\delta \Pi^{cons} 
\Bigl( 
\bigl( \boldsymbol{u}; p^{(f)} \bigr) 
;  
\delta \eta 
\Bigr) 
\biggr) 
; 
\bigl( \delta \boldsymbol{u}, \delta \boldsymbol{p}^{(f)} \bigr) 
\Biggr).
\end{align}
The second variations appearing in the above equations can expressed as follows
\begin{multline}
\delta 
\Biggl(  
\biggl( 
\delta \Pi^{eq} 
\Bigl( 
\bigl( \boldsymbol{u}; p^{(f)} \bigr) 
;  \delta \boldsymbol{w} 
\Bigr) 
\biggr) 
; 
\bigl( \delta \boldsymbol{u}, \delta \boldsymbol{p}^{(f)} \bigr) 
\Biggr)
-
\delta 
\Biggl(
\biggl( 
\delta \Pi^{cons} 
\Bigl( 
\bigl( \boldsymbol{u}; p^{(f)} \bigr) 
;  
\delta \eta 
\Bigr) 
\biggr) 
; 
\bigl( \delta \boldsymbol{u}, \delta \boldsymbol{p}^{(f)} \bigr) 
\Biggr)
\\
=
\begin{bmatrix}
\delta \boldsymbol{w}, \delta \eta 
\end{bmatrix}
\underbrace{\begin{bmatrix}
\mathbf{A} & \mathbf{B}^{(1)}
\\
\mathbf{B}^{(2)} & \mathbf{C}
\end{bmatrix}
}_{\coloneqq \mathbf{J}}
\begin{bmatrix}
\delta \boldsymbol{u}
\\
\delta p^{(f)}
\end{bmatrix},
\end{multline}
where $\mathbf{J}$ indicates the Jacobian matrix. Under the assumptions made so far, the submatrix $\mathbf{C}$ is symmetric, i.e., $\mathbf{C} = \mathbf{C}^{T}.$

In the small strain regime, the above linear system satisfies also the useful properties $\mathbf{A} = \mathbf{A}^T$ and $\left( \mathbf{B}^{(1)} \right)^T = \mathbf{B}^{(2)}$. 
In this context, the Jacobian matrix can exhibit two situations which are well-acknowledged for saddle-point problems. 
On the one hand, in the case of drained or partially drained processes, the above Jacobian is as follows
\begin{equation}
\mathbf{J} = 
\begin{bmatrix}
\mathbf{A} & \mathbf{B}^T
\\
\mathbf{B} & \mathbf{C}
\end{bmatrix}.
\end{equation}
The solvability (uniqueness of solution) of the above linear system is guaranteed if the bilinear forms associated with the submatrices $\mathbf{A}$ and $\mathbf{B}$ are continuous, and $\mathbf{A}$ and $\mathbf{C}$ are coercive (see, for instance, Boffi \emph{et al.}~\cite{boffi2013mixed}, Proposition 4.3.1).
On the other hand, the process can take place in nearly-undrained conditions.
Physically speaking, these circumstances can occur when load rates are rapidly applied with low values of hydraulic conductivity or when BCs do not allow for drainage. 
Theses cases result in the following Jacobian
\begin{equation}
\label{eq:Jacobian_undrained}
\mathbf{J}
=
\begin{bmatrix}
\mathbf{A} & \mathbf{B}^T
\\
\mathbf{B} & \epsilon \, \tilde{\mathbf{C}}
\end{bmatrix},
\end{equation}
where $\mathbf{C} = \epsilon \, \tilde{\mathbf{C}}$ with $0 < \epsilon << 1$ indicates that the submatrix $\mathbf{C}$ contains small entries compared to the other submatrices. 
In the above condition, solvability is guaranteed (see, e.g., Boffi \emph{et al.}~\cite{boffi2013mixed}, Theorem 4.3.4) if the bilinear form associated with $\mathbf{A}$ and $\mathbf{C}$ are continuous and $\mathbf{A}$ is coercive. Furthermore, the solution to Eq.~\eqref{eq:Jacobian_undrained} is stable  if $\mathbf{B}^T$ is such that its bilinear form satisfies the \emph{inf-sup} condition, i.e.,
\begin{equation}
\label{eq:inf-sup}
\inf_{ \left( {\delta \eta} \in {\mathscr{W}_{0}} \right)} \ \sup_{\left(  {\delta \boldsymbol{w}} \in 
{\mathscr{V}_{\boldsymbol{0}}} \right)
} \dfrac{
\int_{\omega} {\delta \eta} \ \text{\textbf{div}} \cdot {\delta \boldsymbol{w}} \, dv 
}{
|| {\delta \eta } ||_{ H^0 }
\ ||{\delta \boldsymbol{w} } ||_{ {\mathscr{V}_{\boldsymbol{0}}} }
}
\geq \alpha > 0,
\end{equation}
where the L-2 norm $|| \left( \bullet \right) ||_{H^0}$ must be considered even for poro-mechanical problems (see the discussions in Mira \emph{et al.}~\cite{mira2003new} and Dortdivanlioglu \emph{et al.}~\cite{dortdivanlioglu2018mixed}). 
Furthermore, the use of this norm in Eq.~\eqref{eq:inf-sup} allows use of numerical test proposed by Chapelle and Bathe~\cite{chapelle1993inf} for incompressible elasticity or Stokes problems.

From the above discussions on drained and undrained cases, it emerges that three requirements are necessary for inverting the Jacobian matrix without incurring in unstable solutions (neglecting for the moment the continuity of bilinear forms of submatrices $\mathbf{A}$ and $\mathbf{C}$). 
These conditions are the coercivity of the bilinear forms associated with $\mathbf{A}$ and $\mathbf{C}$, and the satisfaction of the inf-sup condition for the bilinear form associated with $\mathbf{B}$. 

It must also be stressed that adopting the improved Hencky material defined by Eq.~\eqref{eq:strain_energy} has important consequences for the coercivity of $\mathbf{C}$, since it avoids negative values of the Eulerian porosity, as detailed in Pretti \emph{et al.}~\cite{pretti2024preserving}.
While this criterion does not come from the numerical requirements of the problem under investigation, its physically based nature substantially affects numerical outcomes.
This constraint and its relative consequences hold even in the case of compressible material. In this case, the reader is invited to refer to the work of Nedjar~\cite{nedjar2013formulation,nedjar2014finite} on circumventing negative values in Eulerian porosity.
Owing to the adoption of the improved Hencky material in this work, this physical constraint will not be discussed further.

\section{Stable MPM discretisation}
\label{sec:MPM_discretisation}

This section introduces an MPM $\boldsymbol{u}-p^{(f)}$ formulation and discusses when the conditions for solvability discussed in the previous section are not fulfilled. 
Remedies to mitigate the effects of these losses are also proposed.

\subsection{MPM algorithm}
\label{subsec:MPM_algorithm}
\begin{figure}
\centering
\includegraphics[width=1\textwidth, trim=2.5cm 0cm 0cm 0cm]{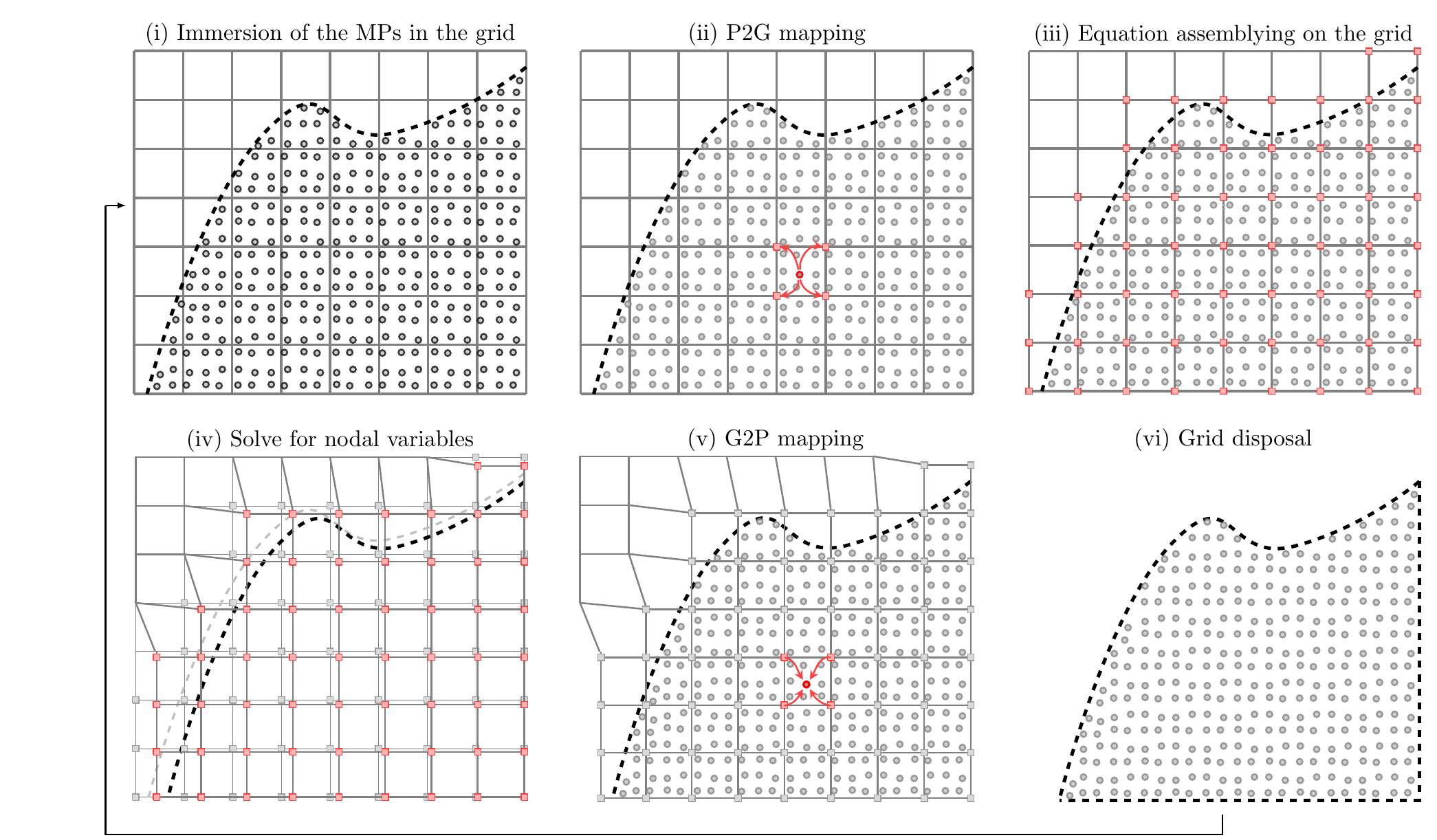}
\caption{MPM step phases.}
\label{fig:MPM_phases}
\end{figure}
Figure~\ref{fig:MPM_phases} shows the phases in the overall MPM algorithm. 
The MPM-based discretisation stipulates that the considered body of global volume $\omega$ is discretised by a cloud of Material Points (MPs), i.e., $\omega \approx \cup_{mp = 1}^{N^{mp}} \  {^{mp}\omega} \coloneqq {^{MP}\omega}$; the MPs carry all the information necessary to run the analysis. 
These points are immersed into a (usually Cartesian) grid/mesh $\mathcal{T}$ discretising a portion of the Euclidean space $\mathcal{E}$ of $n^{dim}$ dimensions, which fully contains the considered body $\omega$ (Phase $(\mathrm{i})$ in Figure~\ref{fig:MPM_phases}). 

Finite-dimensional test and trial functions are defined on a mesh of generic length $\tilde{h}$. $\tilde{h}$, when employed as a superscript, denotes finite-dimensional quantities. 
The trial and test functions functions are interpolated with the help of nodal shape functions, similar to the FEM. The finite-dimensional spaces of the trial functions are given by 
\begin{multline}
\label{eq:space_displacement_generic}
{^{\tilde{h}}\mathscr{V}_{\bar{\boldsymbol{u}}}} 
 = 
\Bigl\{
{^{\tilde{h}} \boldsymbol{u}} 
\in 
\left[ C^{m} \left( ^{\tilde{h}} \bar{\omega} \right) \right]^{n^{dim}}
\bigm| \,
{^{\tilde{h}} \boldsymbol{u}} \left( \boldsymbol{x} \right) = \mathrm{N}^{\tilde{h},u}_{a} \left( \boldsymbol{x} \right) \mathbf{u}_{a} \ \forall a \in \text{clos} \left( {^{\tilde{h}}\mathcal{T}^{act}} \right), \, 
\\
{^{\tilde{h}}\boldsymbol{u}} = \bar{\boldsymbol{u}} 
\ \hbox{on} \ {^{\tilde{h}}\bar{\gamma}^{\bar{u}}} 
\Bigr\};
\end{multline}
\begin{equation}
 \label{eq:space_pressures_generic}
{^{\tilde{h}} \mathscr{W}_{\bar{p}}} 
 =
\left\{
^{\tilde{h}}p^{(f)} \in  C^{m} \left( ^{\tilde{h}}\bar{\omega} \right)
\bigm| \,
{^{\tilde{h}} p^{(f)}} \left( \boldsymbol{x} \right) = \mathrm{N}^{\tilde{h}, p}_{a} \left( \boldsymbol{x} \right) \,  \mathrm{p}_{a} \ \forall a \in \text{clos} \left( {^{\tilde{h}}\mathcal{T}^{act}} \right)
\right\},
\end{equation}
where $^{\tilde{h}}\bar{\omega} = \left( \bigcup_{T \in {^{\tilde{h}}\mathcal{T}^{act}} } T \right)$ denotes the volume defined by the set of active grid elements ${^{\tilde{h}}\mathcal{T}^{act}}$, i.e.,
\begin{equation}
{^{\tilde{h}}\mathcal{T}^{act}} 
= 
\bigl\{ T \in {^{\tilde{h}}\mathcal{T}} \bigm|  \exists \, mp:  \mathrm{N}^{\tilde{h}}_a \left( ^{mp}\boldsymbol{x}\right) > 0 \ \text{with} \, a \in \text{clos} (T) \bigr\} ,
\end{equation}
with $\mathrm{N}^{\tilde{h}}_a$ being the generic low-order shape functions at the $a-$th node belonging to the element $T$ and $\text{clos} \left( \bullet \right)$ denotes the closure of the entity $\left( \bullet \right)$\footnote{
This work often refers to the closure of a generic entity, which is a shorthand for indicating the set of entities contained by the considered generic one and the entity itself. If, for instance, a bi-dimensional mesh element $T$ is considered, the closure of $T$ includes the faces and nodes belonging to $T$ and $T$ itself.
}. 
As detailed below in Section~\ref{subsec:Q1-iso-Q2-Q1_discretisation}, this work assumes two kinds of shape functions. 
The first choice is given by first-order Lagrange polynomials, and results in $C^{0} ( {^{\tilde{h}} \bar{\omega}} )$ functions (i.e., $m = 0$). Since they were first adopted in the original MPM formulation (see Sulsky \emph{et al.}~\cite{sulsky1994particle,sulsky1995application}), these shape functions are labelled sMPM (standard MPM). 
On the other hand, GIMPM shape functions (Generalised Interpolation MPM, see, for instance,~\cite{bardenhagen2004generalized,charlton2017igimp}), resulting in piece-wise first-order and second-order polynomials, are also considered, giving $C^{1} ( {^{\tilde{h}} \bar{\omega}} )$ functions (i.e., $m = 1$).  

To initialise the grid with the required information, a mapping from the MPs to the grid nodes is usually required (Phase $(\mathrm{ii})$ in Figure~\ref{fig:MPM_phases}). 
This procedure is named Point-2-Grid (P2G) mapping and, for the current formulation, is addressed in Section~\ref{subsec:mapping_processes}.

The (standard) assembly process of the Lagrangian equations (Phase $(\mathrm{iii})$ in Figure~\ref{fig:MPM_phases}) takes place at the grid nodes and employs the MPs as quadrature points. 
In the case of the considered $\boldsymbol{u}-p^{(f)}$ formulation, the equations become (dropping the dependencies from the position) 
\begin{multline}
\label{eq:eq_equation_discrete_form}
\int_{{^{MP}\omega}} \updelta \mathbf{w}_{a} \,\text{\textbf{grad}} \left( \mathrm{N}^{\tilde{h}, u}_{a} \right) \boldsymbol{:} \left( \boldsymbol{\sigma}'
- \mathrm{N}^{\tilde{h},p}_{b} \, \mathrm{p}_{b} \, \boldsymbol{1} \right) \ d v 
-
\int_{{^{MP}\omega}}  \left( \updelta \mathbf{w}_{a} \,  \mathrm{N}^{\tilde{h}, u}_{a} \right) \cdot \rho \boldsymbol{b} \ dv 
\\
-
\int_{\gamma^{\bar{t}}}  \left( \updelta \mathbf{w}_{a} \,  \mathrm{N}^{\tilde{h},u}_{a} \left( \boldsymbol{x} \right) \right) \cdot \bar{\boldsymbol{t}} \ da
= 0, 
\quad \forall \updelta \mathbf{w}_{a} \in \ ^{\tilde{h}}\mathscr{V}_{\boldsymbol{0}};
\end{multline}
\begin{multline}
\label{eq:fluid_mass_conservation_discrete_form}
\int_{{^{MP}\omega}} \left( \updelta \upeta_{a} \,  \mathrm{N}^{\tilde{h},p}_{a}  \right) \rho_0^{(f)}  \frac{\dot{J}}{J} \ dv 
+
\int_{{^{MP}\omega}} \frac{\upkappa}{g} \, \updelta \upeta_{a} \, \text{\textbf{grad}}  \left(  \mathrm{N}^{\tilde{h},p}_{a}\right) \cdot  \left( \text{\textbf{grad}} \left( \mathrm{N}^{\tilde{h},p}_{b} \right) \mathrm{p}_{b} - \rho^{(f)} \boldsymbol{b} \right) dv
\\
-
\int_{\gamma^{\bar{q}}} \left( \updelta \upeta_{a} \,  \mathrm{N}^{\tilde{h},p}_{a}   \right) \, \bar{q} \ da 
+ 
\gamma^{pen}
\int_{\gamma^{\bar{p}}} \left( \updelta \upeta_{a} \,  \mathrm{N}^{\tilde{h},p}_{a}   \right) \, \left(  \mathrm{N}^{\tilde{h},p}_{b} \, \mathrm{p}_{b} - \bar{p} \right) da
= 0,
\quad \forall \updelta \upeta_{a} \in \ {^{\tilde{h}}\mathscr{W}_{0}},
\end{multline}
where the volume integrals are approximated as follows 
\begin{equation}
\label{eq:numerical_integration}
\int_{{^{MP}\omega}} \left( \bullet \right) dv \approx
\sum_{mp = 1}^{N^{mp}} {^{mp}\omega} \left( \bullet \left( ^{mp}\boldsymbol{x} \right) \right).
\end{equation}
As for the boundary integrals, this work follows the procedure explained by Bird \emph{et al.}~\cite{bird2024implicit}, leveraging the MPs' corners (see Bird \emph{et al.}~\cite{bird2024implicit} for more details).
This method is also exploited for non-conforming pressure Dirichlet BCs, represented by the penalty term in Eq.~\eqref{eq:fluid_mass_conservation_discrete_form}, with $\gamma^{pen}$ being a penalty parameter.
Conversely, this work always considers mesh-conforming displacement BCs, as assumed in Eq.~\eqref{eq:space_displacement_generic}.
Eqs.~\eqref{eq:eq_equation_discrete_form} and~\eqref{eq:fluid_mass_conservation_discrete_form} are solved for the nodal unknowns (Phase $(\mathrm{iv})$ in Figure~\ref{fig:MPM_phases}), and this work considers an interactive Newton-Raphson (NR) scheme based on the discrete counterpart of Eq.~\eqref{eq:linearised_functional}.
Since the updated solution lies on the grid, a subsequent mapping from this to the MPs is required to update the solution at the MPs level (Phase $(\mathrm{v})$ in Figure~\ref{fig:MPM_phases}). 
This is called the Grid-2-Point (G2P) mapping and is discussed in Section~\ref{subsec:mapping_processes}.
The grid is discarded at the end of the step (Phase $(\mathrm{vi})$ in Figure~\ref{fig:MPM_phases}), and a further step can start with an undistorted grid.

\subsection{The small-cut issue in the MPM}
\label{subsec:small-cut_issue}
As will be clear from the description of the phases of the MPM, the immersion of the clouds of MPs into a mesh and their use as quadrature points generates an intrinsically unfitted method since ${^{MP}\omega} \subseteq {^{\tilde{h}}\bar{\omega}}$. 
Moreover, an extremely small overlap between the shape function's stencil and the MP's domains can lead to ill-conditioned systems.
This problem, named the small-cut issue is exacerbated especially for those elements on the domain's boundaries.
As demonstrated by Burman~\cite{burman2010ghost} in the context of unfitted FEM, this issue provokes the loss of the coercivity for the bilinear forms computed on these elements and the consequent impossibility of solving the relative linear (or linearised) system of equations. 

To manage this issue, Burman~\cite{burman2010ghost} has proposed adding a penalty term, named the ghost penalty stabilisation, which extends the coercivity of the bilinear forms to the partially filled boundary elements . 
While strictly enforcing coercivity in an MPM formulation is difficult to prove (properties such as continuity or coercivity rely on norms which, in turn, are based on quadrature rules more precise than Eq.~\eqref{eq:numerical_integration}), it is nonetheless possible to limit the condition number of the matrices by the use of the ghost penalty method (see Coombs~\cite{coombs2023ghost}). 
The addition of the ghost penalty restores to some degree coercivity on those elements suffering from the small-cut issue.
Hence, while formalisms must be dropped for the above-mentioned causes, MPM-based discrete linear systems can be inverted independently from the small-cut issue, as long as the ghost penalty (or a similar technique, see the discussion in Section~\ref{subsec:face_ghost_stabilisations}) is introduced.

Among the ghost penalties introduced in the literature for the unfitted FEM (see Burman \emph{et al.}~\cite{burman2022design}), the \emph{face} ghost penalty\footnote{
While the \emph{face} ghost penalty is the term usually adopted in the literature, this work refers to \emph{facet} to indicate, in the $n^{dim}-$th dimensional space, the face with dimension $n^{dim}-1$, i.e., a point, an edge, and a face for $n^{dim} = 1, 2, 3$, respectively.
} has been introduced to the MPM by Coombs~\cite{coombs2023ghost}, and an extension to a $\boldsymbol{u}-p^{(f)}$ formulation for the unfitted FEM has been proposed by Liu \emph{et al.}~\cite{liu2022unfitted}.
The following section extends the face ghost stabilisation for the $\boldsymbol{u}-p^{(f)}$ MPM formulation introduced above.

\subsection{Face ghost stabilisation}
\label{subsec:face_ghost_stabilisations}
In the case of low-order shape functions, the face ghost stabilisation states that the following bilinear form must be added to the considered primary equations, i.e.,
\begin{equation}
\label{eq:face_ghost}
j^{F} \left( {^{\tilde{h}}\delta \boldsymbol{v}}, {^{\tilde{h}}\boldsymbol{v}} \right) 
\coloneqq
\tilde{h}
\int_{^{\tilde{h}}\gamma^{F}} 
\left[ \left[
\left( \text{\textbf{grad}} \left( {^{\tilde{h}}\delta \boldsymbol{v}}  \right) \cdot \mathbf{n} \right) 
\right] \right]
\cdot 
\left[ \left[
\left( \text{\textbf{grad}} \left(  ^{\tilde{h}}\boldsymbol{v}  \right) \cdot \mathbf{n} \right) 
\right] \right] \ da,
\end{equation}
where $^{\tilde{h}}\boldsymbol{v}$ and $^{\tilde{h}}\delta \boldsymbol{v}$ are the generic trial and test functions and $[[ \text{\textbf{grad}}\left( \bullet \right) \, \mathbf{n}]]$ is the jump in the normal gradient across a shared facet $F$. 
The procedure to select these facets (listed in $\mathcal{F}$) describing the surface $^{\tilde{h}}\gamma^F$ is outlined in Algorithms~\ref{algorithm:boundary_els_selection} and~\ref{algorithm:faces_selection}.
In particular, since the MPM retains no explicit representation of the boundary, Algorithm~\ref{algorithm:boundary_els_selection} is necessary to track which grid elements constitute the domain's boundary ${^{\tilde{h}}\mathcal{T}^{\partial {^{\tilde{h}} }\bar{\omega}}}$. 
To read this algorithm, it is necessary also to introduce the set of inactive elements $^{\tilde{h}}\mathcal{T}^{int}$, defined as $^{\tilde{h}}\mathcal{T}^{in} \coloneqq {^{\tilde{h}}\mathcal{T}} \setminus {^{\tilde{h}}\mathcal{T}^{act}}$.
Among the facets belonging to the set of boundary elements ${^{\tilde{h}}\mathcal{T}^{\partial {^{\tilde{h}} }\bar{\omega}}}$, Algorithm~\ref{algorithm:faces_selection} selects those satisfying one of the following criteria: 
\begin{itemize}
    \item the facet is shared between one boundary element and one active non-boundary element; or
    \item the facet is between two adjacent boundary elements.
\end{itemize}
While selecting the former set of facets (criterion at the first bullet point) is justified by prolonging the coercivity from the bulk of the material, the latter (criterion at the second bullet point) is mainly required for boundary corner elements not directly connected to an active non-boundary element.
The union of the selected facets, i.e.,
\begin{equation}
^{\tilde{h}}\gamma^{F} =  \bigcup_{F \in \mathcal{F}} F 
\end{equation}
constitutes the surface $^{\tilde{h}}\gamma^F$ necessary for Eq.~\eqref{eq:face_ghost}. Unlike the  volume integrals in Eq.~\eqref{eq:numerical_integration}, Eq.~\eqref{eq:face_ghost} is integrated via Gauss-Legendre quadrature rule using Gauss Points seeded on the grid.

\begin{algorithm}[]
\footnotesize{
${^{\tilde{h}}\mathcal{T}^{\partial {^{\tilde{h}} }\bar{\omega}}} \leftarrow []$\;
el\_count $\leftarrow$ 1\;
 \ForEach(
 \tcp*[f]{\scriptsize{loop over active elements}}
 ){
 element $T \in {^{\tilde{h}}\mathcal{T}^{act}}$ 
 }{ 
 face\_count $\leftarrow$ 1\;
    \ForEach(
    \tcp*[f]{\scriptsize{loop over facets of the selected active element}}
    ){
    facet $F \in \text{clos} (T)$
    }{
    \If(
    \tcp*[f]{\scriptsize{check if facet is also shared with an inactive element}}
    ){
    $F \in \text{clos} ({^{\tilde{h}}\mathcal{T}^{in}})$ 
    }{
    Add $T$ to ${^{\tilde{h}}\mathcal{T}^{\partial {^{\tilde{h}} }\bar{\omega}}}$;
    }
     face\_count ++
    }
    el\_count ++
 }
}
\caption{Selection of the boundary elements belonging to ${^{\tilde{h}}\mathcal{T}^{\partial {^{\tilde{h}} }\bar{\omega}}}$.}
\label{algorithm:boundary_els_selection}
\end{algorithm}

\begin{algorithm}[]
\footnotesize{
$\mathcal{F} \leftarrow []$\;
el\_count $\leftarrow$ 1\;
 \ForEach(
 \tcp*[f]{\scriptsize{loop over boundary elements}}
 ){
 element $T_i \in {^{\tilde{h}}\mathcal{T}^{\partial {^{\tilde{h}} }\bar{\omega}}}$
 }{
   face\_count $\leftarrow$ 1\;
    \ForEach(
    \tcp*[f]{\scriptsize{loop over facets of the selected boundary element}}
    ){
    facet $F \in \text{clos} ( T_i )$
    }{
    \If(
    \tcp*[f]{\scriptsize{check if facet is shared with another boundary element}}
    ){
    $F \in \text{clos} \left( {^{\tilde{h}}\mathcal{T}^{\partial {^{\tilde{h}} }\bar{\omega}}}  \setminus T_i \right)$
    }{
    Add $F$ to $\mathcal{F}$;
    }
    \If(
    \tcp*[f]{\scriptsize{check if facet is shared with active non-boundary element}}
    ){
    $F \in \text{clos} \left( {^{\tilde{h}}\mathcal{T}^{act} }  \setminus {^{\tilde{h}}\mathcal{T}^{\partial {^{\tilde{h}} }\bar{\omega}} }  \right)$
    }{
    Add $F$ to $\mathcal{F}$;
    }
    face\_count ++
    }
 el\_count ++
 }
}
\caption{Selection of the facets for the face ghost penalty.}
\label{algorithm:faces_selection}
\end{algorithm}

Algorithms~\ref{algorithm:boundary_els_selection} and~\ref{algorithm:faces_selection} (as well as Eq.~\eqref{eq:face_ghost}) are relevant to sMPM and GIMPM shape functions. 
Unlike the unfitted FEM, in which the boundary elements constitute a single layer of elements, high-order MPMs can contain multiple layers of boundary elements due to the stencil of these high-order shape functions, which cover multiple elements simultaneously. This situation can also occur for the GIMPM.
However, since only two layers of boundary elements can be active simultaneously in the worst-case scenario for the GIMPM, Algorithms~\ref{algorithm:boundary_els_selection} and~\ref{algorithm:faces_selection} prove to be sufficient (see Example~\ref{subsec:investigation_ghost}). 

As stabilising multiple layers can be quite burdensome, Yamaguchi \emph{et al.}~\cite{yamaguchi2021extended} adapted the Extended B-splines (EBS, see H{\"o}llig \emph{et al.}~\cite{hollig2001weighted,hollig2003finite}) technique to stabilise MPM formulations employing B-spline shape functions. 
Another interesting technique has been proposed by Badia \emph{et al.}~\cite{badia2018aggregated} for the unfitted FEM, which aggregates elements on the boundary suffering from the small-cut issue.  
However, these techniques can be problematic to integrate into the current framework: on the one hand, the EBS-MPM cannot be used for the considered low-order shape functions, while aggregating can be particularly difficult, especially in light of the inf-sup stable overlapping meshes adopted in this manuscript. These are explained in the following section.



\subsection{Qk$_{SD}$-Qk spatial discretisation}
\label{subsec:Q1-iso-Q2-Q1_discretisation}
To comply with the inf-sup requirement, the concept of Qk$_{SD}$-Qk elements (where the subscript $SD$ stands for \emph{S}ub-\emph{D}ivided, while $k$ is the degree of the shape functions polynomials) from the FEM is here adapted to the MPM. 
This is achieved by overlapping a coarser grid for the pressure field (of length $h$) on a finer mesh for the displacement field (of length $h/2$). For simplicity, both these meshes are considered as Cartesian. 

Qk$_{SD}$-Qk elements are well-acknowledged in the FEM (see, for instance, Dortdivanlioglu \emph{et al.}~\cite{dortdivanlioglu2018mixed} for a $\boldsymbol{u}-p^{f}$ formulation), and the idea can be traced at least back to Bercovier and Pironneau~\cite{bercovier1979error}, who gave the first error-estimate for two- and three-dimensional FEM of these kind of elements. 
In the case of a MPM formulation, to the best of authors' knowledge, Qk$_{SD}$-Qk concept has been applied only to problems involving incompressible elasticity and finite elements with B-spline shape functions ~\cite{madadi2024subdivision}.

This choice implies that the finite-dimensional spaces~\eqref{eq:space_displacement_generic} and~\eqref{eq:space_pressures_generic} for the trial functions become
\begin{multline}
{{^{h/2}}\mathscr{V}_{\bar{\boldsymbol{u}}}} 
= 
\Bigl\{
{^{h/2} \boldsymbol{u}} 
\in 
\left[ C^{m} \left( ^{h/2} \bar{\omega} \right) \right]^{n^{dim}}
\bigm| \,
{^{h/2} \boldsymbol{u}} \left( \boldsymbol{x} \right) = \mathrm{N}^{h/2}_{a} \left( \boldsymbol{x} \right) \mathbf{u}_{a} 
\\
\forall a \in \text{clos} \left( {^{h/2}\mathcal{T}^{act}} \right), \,
{^{h/2}\boldsymbol{u}} = \bar{\boldsymbol{u}} 
\ \hbox{on} \ {^{h/2}\bar{\gamma}^{\bar{u}}} 
\Bigr\};
\end{multline}
\begin{equation}
{{^h}\mathscr{W}_{\bar{p}}} 
=
\left\{
^{h}p^{(f)} \in  C^{m} \left( ^{h}\bar{\omega} \right)
\bigm| \,
{^{h} p^{(f)}} \left( \boldsymbol{x} \right) = \mathrm{N}^{h}_{a} \left( \boldsymbol{x} \right) \,  \mathrm{p}_{a} \ \forall a \in \text{clos} \left( {^{h}\mathcal{T}^{act}} \right)
\right\},
\end{equation}
and similar definitions follow for the spaces of test functions ${^{h/2}}\mathscr{V}_{\boldsymbol{0}}$ and ${^h}\mathscr{W}_{\boldsymbol{0}}$.
Using Qk$_{SD}$-Qk-inspired meshes and including the face ghost stabilisation Eq.~\eqref{eq:face_ghost}, the weak form becomes: find ${^{h/2}\boldsymbol{u}} \in {^{h/2}\mathscr{V}_{\bar{\boldsymbol{u}}}}$ and ${^{h}p^{(f)}} \in {{^h}\mathscr{W}_{\bar{p}}}$ such that, for $t \in [ 0, \tilde{t} ] $,
\begin{multline}
\label{eq:eq_equation_discrete_form_QaSD-Qa}
\int_{{^{MP}\omega}} \updelta \mathbf{w}_{a} \,\text{\textbf{grad}} \left( \mathrm{N}^{h/2}_{a} \right) \boldsymbol{:} \left( \boldsymbol{\sigma}'
- \mathrm{N}^{h}_{b} \, \mathrm{p}_{b} \, \boldsymbol{1} \right) \ d v 
\\
-
\int_{{^{MP}\omega}}  \left( \updelta \mathbf{w}_{a} \,  \mathrm{N}^{h/2}_{a} \right) \cdot \rho \boldsymbol{b} \ dv 
-
\int_{\gamma^{\bar{t}}}  \left( \updelta \mathbf{w}_{a} \,  \mathrm{N}^{h/2}_{a} \right) \cdot \bar{\boldsymbol{t}} \ da
\\ 
+
\gamma^A \, 
\frac{h}{2}
\int_{{^{h/2}}\gamma^{F}} \updelta \mathbf{w}_{a} \left[ \left[ \left( \text{\textbf{grad}} \left( \mathrm{N}^{h/2}_{a}  \right) \cdot \mathbf{n} \right) \right] \right]
\cdot 
\left[ \left[ \left( \text{\textbf{grad}} \left( \mathrm{N}^{h/2}_{b} \right) \cdot \mathbf{n}   \right ) \right] \right] \mathbf{u}_{b} \ da
= 0, 
\\
\quad \forall \updelta \mathbf{w}_{a} \in \ ^{h/2}\mathscr{V}_{\boldsymbol{0}};
\end{multline}
\begin{multline}
\label{eq:fluid_mass_conservation_discrete_form_QaSD-Qa}
\int_{{^{MP}\omega}} \left( \updelta \upeta_{a} \,  \mathrm{N}^{h}_{a}  \right) \rho_0^{(f)}  \frac{\dot{J}}{J} \ dv 
+
\int_{{^{MP}\omega}} \frac{\upkappa}{g} \, \updelta \upeta_{a} \, \text{\textbf{grad}}  \left(  \mathrm{N}^{h}_{a}\right) \cdot  \left( \text{\textbf{grad}} \left( \mathrm{N}^{h}_{b} \right) \mathrm{p}_{b}  - \rho^{(f)} \boldsymbol{b} \right) dv
\\
-
\int_{\gamma^{\bar{q}}} \left( \updelta \upeta_{a} \,  \mathrm{N}^{h}_{a}   \right) \, \bar{q} \ da
+
\gamma^{pen}
\int_{\gamma^{\bar{p}}} \left( \updelta \upeta_{a} \,  \mathrm{N}^{h}_{a}   \right) \, \left(  \mathrm{N}^{h}_{b} \, \mathrm{p}_{b} - \bar{p} \right) da
\\
+
\,
\gamma^C 
h
\int_{^{h}\gamma^{F}} \updelta \upeta_{a} 
\left[ \left[ 
\left( \text{\textbf{grad}} \left( \mathrm{N}^{h}_{a}  \right) \cdot \mathbf{n} \right) 
\right] \right]
\left[ \left[ 
\left( \text{\textbf{grad}} \left( \mathrm{N}^{h}_{b} \right) \cdot \mathbf{n} \right) 
\right] \right] \mathrm{p}_{b} 
\ da
= 0,
\quad \forall \updelta \upeta_{a} \in \ {{^h}\mathscr{W}_{0}},
\end{multline}
giving the initial condition ${^{h/2}\boldsymbol{u}} \left( t = 0 \right) = \boldsymbol{u}_0$.
In the above equations, $\gamma^A$ and $\gamma^C$ are the coefficients used to the scale the effect of the ghost penalty terms and to match the physical units of the equations. For this latter purpose, the units of $\gamma^A$ must be similar to a stress measure, while the units of $\gamma^C$ should be the same of the $\upkappa/g$.

\begin{figure}[!h]
\centering
\begin{subfigure}[t]{0.495\textwidth}
\includegraphics[width=\textwidth,trim= 1.2cm 0.15cm 0 0]{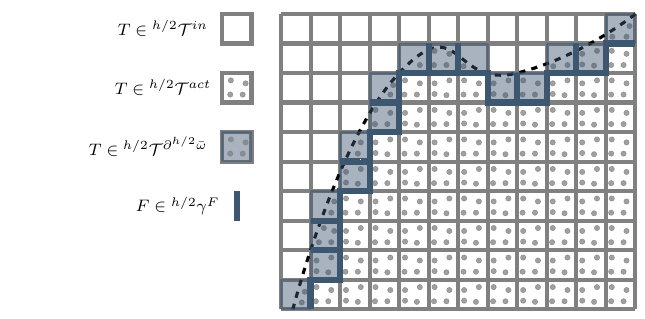}
\caption{Application of Algorithms~\ref{algorithm:boundary_els_selection} and~\ref{algorithm:faces_selection} to the finer mesh for $^{h/2}\boldsymbol{u}$.}
\label{subfig:finer_mesh}
\end{subfigure}
\hfill
\begin{subfigure}[t]{0.495\textwidth}
\includegraphics[width=\textwidth,trim= 1.2cm 0.15cm 0 0]{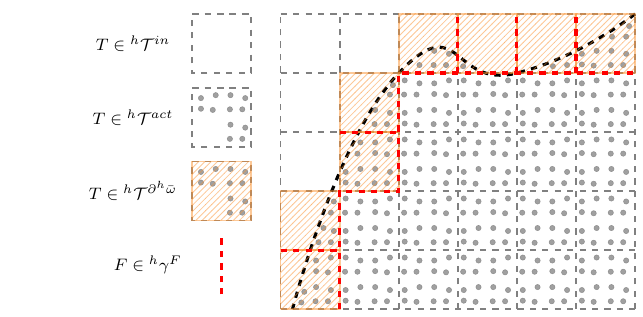}
\caption{Application of Algorithms~\ref{algorithm:boundary_els_selection} and~\ref{algorithm:faces_selection} to the coarser mesh for $^hp^{(f)}$.}
\label{subfig:coarser_mesh}
\end{subfigure}
\\
\begin{subfigure}[t]{0.35\textwidth}
\includegraphics[width=\textwidth,trim= 1.2cm 0.15cm 0 0]{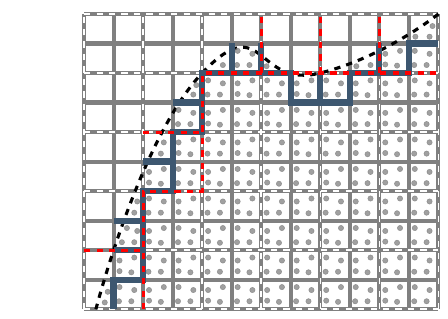}
\caption{Final result of facets where face ghost stabilisations are applied.}
\label{subfig:facets_stabilised}
\end{subfigure}
\caption{Algorithmic selection of boundary elements and facets (top row) for finer (left) and coarser (right) meshes, and stabilised facets (bottom row).}
\label{fig:ghosts}
\end{figure}

Fig.~\ref{subfig:finer_mesh} shows the application of the face ghost stabilisation to the finer mesh of length $h/2$ for the displacement field (i.e., ghost penalty term added to Eq.~\eqref{eq:eq_equation_discrete_form_QaSD-Qa}), distinguishing between inactive and active elements, boundary elements and boundary facets. 
Similarly, Fig.~\ref{subfig:coarser_mesh} illustrates the same quantities for the coarser mesh of length $h$ employed for the pressure field (i.e., ghost penalty term added to Eq.~\eqref{eq:fluid_mass_conservation_discrete_form_QaSD-Qa}). 
The combined result of the different stabilised surfaces for the displacement and pressure fields are represented in Fig.~\ref{subfig:facets_stabilised}.

For the sake of completeness, it must be acknowledged that a plethora of inf-sup stabilised elements exist in the literature and covering the complete literature on this subject is beyond the scope of this work. 
As for the FEM, the reader can refer to, for instance, Leborgne~\cite{leborgne2023inf} for a recent publication discussing the inf-sup condition for different mixed formulations (and possible remedies) or to the monograph of Boffi \emph{et al.}~\cite{boffi2013mixed} for a complete discussion. 
The nature of Qk$_{SD}$-Qk elements makes its adaptation straightforward to use in the MPM, especially if low-order shape functions are used. 
On top of this, stable elements by design possess the advantage of not adding terms to the primary equations. 
This is particularly advantageous if the NR method is employed to solve the implicit primary equations, since less terms can be considered in the Jacobian matrix. 
In the following section, Qk$_{SD}$-Qk meshes are diversified for the sMPM and the GIMPM and their inf-sup stability is investigated.

\subsubsection{sMPM shape functions}
\label{subsubsec:sMPM}
Having selected a Cartesian mesh composed by quadrilateral (or hexahedral, if $n^{dim} = 3$) elements, the basis functions are constructed using the tensor product of one-dimensional functions. 
In the case of sMPM for the generic mesh of length $\tilde{h}$, first-order Lagrange polynomials are employed, i.e.,
\begin{alignat}{2}
\label{eq:linear_MPM_SFs}
\mathrm{N}^{\tilde{h},1} \left( \xi \right) & = 
\begin{cases}
1 + \xi/\tilde{h}, 
\qquad & 
\hbox{if} \quad -\tilde{h} < \xi \leq 0; 
\\
1 - \xi/\tilde{h}, 
\qquad &
\hbox{if} \quad 0 < \xi \leq \tilde{h},
\end{cases}
\end{alignat}
where $\xi = {^{mp}x} - {^{\tilde{h}}x_v}$ denotes the one-dimensional local coordinate computed as the difference between the material point's coordinate and the grid vertex coordinate. 
The use of these shape functions for a Q1$_{SD}$-Q1 FEM $\boldsymbol{u}-p^{(f)}$ formulation has been tested in Dortdivanlioglu \emph{et al.}~\cite{dortdivanlioglu2018mixed}, employing the patch test proposed by Chapelle and Bathe~\cite{chapelle1993inf} for mixed FEM formulations.
Employing this result and given the similarities between the FEM and the MPM, it can be concluded that Q1$_{SD}$-Q1 mesh-based solutions for sMPM discretisation inherit inf-sup stability from its FEM counterpart. 
This conclusion is also substantiated by Example~\ref{subsec:Terzaghi_consolidation}.

\subsubsection{GIMPM shape functions}
\label{subsubsec:GIMPM}
The GIMPM shape functions in 1D are obtained via the convolution of a constant function denoting the length of the particle $2 l_p$ (named the \emph{characteristic function}) with the sMPM shape function given by Eq.~\eqref{eq:linear_MPM_SFs}.
The reader is referred to Bardenhagen \emph{et al.}~\cite{bardenhagen2004generalized} or Charlton \emph{et al.}~\cite{charlton2017igimp} for this convolution integral procedure.
The one-dimensional GIMPM shape functions defined on a mesh of length $\tilde{h}$ is given by
\begin{alignat}{2}
\label{eq:linear_GIMPM_SFs}
\mathrm{S}^{\tilde{h}} \left( \xi \right) & = 
\begin{cases}
\left(\tilde{h} + l_p + \xi \right)^2/\left( 4 \tilde{h} l_p \right),
\qquad & 
\hbox{if} \quad -\tilde{h}-l_p < \xi \leq -\tilde{h} + l_p;  
\\
1 +  \xi/\tilde{h}, 
\qquad & 
\hbox{if} \quad -\tilde{h}+l_p < \xi \leq - l_p; 
\\
1 -  \xi ^2/\left( 2 \tilde{h} l_p \right) - l_p^2/\left( 2 \tilde{h} l_p \right),
\qquad & 
\hbox{if} \quad -l_p < \xi \leq  l_p; 
\\
1 - \xi/\tilde{h},
\qquad & 
\hbox{if} \quad l_p < \xi \leq \tilde{h} - l_p,
\\
\left(\tilde{h} + l_p - \xi \right)^2/\left( 4 \tilde{h} l_p \right),
\qquad & 
\hbox{if} \quad \tilde{h}-l_p < \xi \leq \tilde{h} + l_p,
\end{cases}
\end{alignat}
with $\xi = {^{mp}x} - {^{\tilde{h}}x_v}$ as above.

\begin{figure}
\centering
\includegraphics[width=0.5\textwidth, trim=0cm 2.85cm 0cm 0cm]{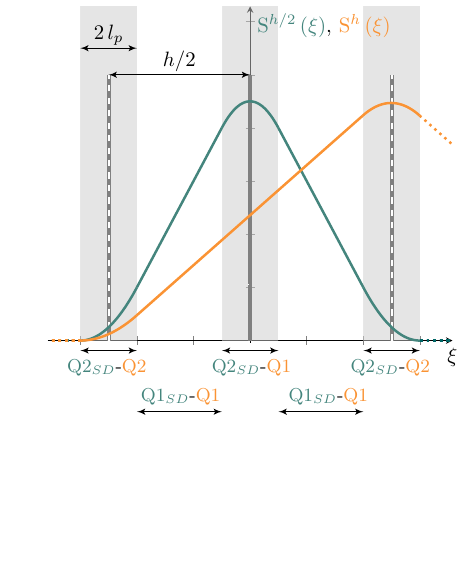}
\caption{Combination of GIMPM functions on overlapping meshes. Thick grey lines illustrate the finer mesh, while white dashed lines represent the coarser one.}
\label{fig:QaSD-Qa_GIMPM}
\end{figure}

Fig.~\ref{fig:QaSD-Qa_GIMPM} shows that three combinations of the polynomials are locally available when using the GIMPM shape functions for Qk$_{SD}$-Qk meshes: Q2$_{SD}$-Q2, Q1$_{SD}$-Q1, and Q2$_{SD}$-Q1. 
As for the case of linear polynomials in Section~\ref{subsubsec:sMPM}, these elements have been separately tested by Dortdivanlioglu \emph{et al.}~\cite{dortdivanlioglu2018mixed} for the FEM using the patch test in Chapelle and Bathe~\cite{chapelle1993inf}.
When examined independently, each of these elements has been shown to be inf-sup stable.
While it is not possible to apply this test directly to an MPM formulation (patch tests cannot be adapted to the MPM), a thought experiment can be built by employing shape functions Eq.~\eqref{eq:linear_GIMPM_SFs} to construct a FEM discretisation and using the results from Dortdivanlioglu \emph{et al.}~\cite{dortdivanlioglu2018mixed}.
The test proposed by Bathe and Chapelle~\cite{chapelle1993inf} checks an eigenvalue-eigenvector problem using matrices relative to discrete bilinear forms. This test is repeated for different decreasing mesh sizes, and it is considered as passed when the minimum eigenvalue does not decrease with finer grids.
The same matrices for the eigenvalue-eigenvector problem employing shape functions in Eq.~\eqref{eq:linear_GIMPM_SFs} are a linear combination of the above elements. This holds since each integration point singularly contributes as Q2$_{SD}$-Q2 or Q1$_{SD}$-Q1 or Q2$_{SD}$-Q1 to its element and the assembling process sums each integration points' contribution. 
Hence, the boundedness of the minimum eigenvalue resulting from this FEM-like discretisation using shape functions in Eq.~\eqref{eq:linear_GIMPM_SFs} follows since it is linear algebra (eigenvalue problem) applied to matrices of bilinear forms.
This rationale is confirmed by numerical Example~\ref{subsec:flexible_strip_footing}, which adopts GIMPM shape functions and exhibits stable pressure values.

\subsection{Time discretisation}
In addition to the spatial discretisation introduced above, a temporal discretisation is required  and this work considers the following Backward-Euler time integration relative to the divergence of the velocity field
\begin{equation}
\label{eq:J_time_discretisation}
\text{\textbf{div}} \cdot \boldsymbol{v}= 
\frac{\dot{J}}{J} 
=
\left( \ln J \right)^{\cdot} 
\approx  
\frac{\left( \ln J_{n+1} - \ln J_n \right)}{\Delta t},
\end{equation}
where the subscript $n+1$ denotes the current time-step and $n$ the previous one.

The literature (see, for instance, Sun \emph{et al.}~\cite{sun2013stabilized} or Zhao and Choo~\cite{zhao2020stabilized}) recognises the role played by the above time-discretisation to avoid negative values of the Jacobian.
However, the current formulation considering an incompressible solid phase imposes a more severe constraint on the Jacobian than $J>0$. 
This can quickly verified if inequalities~\eqref{ineq:bounded_porosity} relative to the Eulerian porosity are substituted in solid mass conservation Eq.~\eqref{eq:Eulerian_porosity}. 
The reader can refer to Pretti \emph{et al.}~\cite{pretti2024preserving} for further details and consequences of the above inequalities.
Despite this consideration, the above formula possesses another desirable feature for the MPM, discussed below in Section~\ref{subsec:mapping_processes} and it is thus employed.

\subsection{Mapping processes and MPs' domain update}
\label{subsec:mapping_processes}
As mentioned in Section~\ref{subsec:MPM_algorithm}, two mapping processes take place in an MPM step. 
The Point-to-Grid (P2G) mapping (phase $(\mathrm{ii})$ in Figure~\ref{fig:MPM_phases}) initialises the nodal grid unknowns, transferring the information from the MPs. 
This mapping process is required to conserve physical properties of interest relative to the mapped information. 
In the context of the dynamics of a solid body, for instance, the velocity mapped from the MPs to the grid is expected to conserve, as much as possible, momenta (linear and angular) and kinetic energy to avoid numerical dissipation (see Love and Sulsky~\cite{love2006unconditionally} or Pretti \emph{et al.}~\cite{pretti2023conservation}). 
A similar rigorous procedure should be expected for a poro-mechanical problem. However, mapping the velocity can be avoided for a quasi-static problem using Eq.~\eqref{eq:J_time_discretisation}, since the Jacobian $J$ is not required on the grid for computational purposes. 
As for the other primary variable, the fluid pressure relates to a form of energy stored in the fluid body only if the considered material is compressible, as pointed out by Miehe \emph{et al.}~\cite{miehe2015minimization}. 
Owing to Assumption~\eqref{hyp:incompressiblity}, this work considers an incompressible fluid constituent, which is not required to be mapped as the pressure acts as a Lagrange multiplier, which does not entail any form of energy. 
Overall, it emerges that G2P mapping is not required in this context, i.e., quasi-static simulations considering an incompressible fluid phase.
However, this procedure can take place for estimating better predictors for the NR scheme, but these estimates are free from numerical dissipation .

Once the grid solution is achieved, this must be mapped to the MPs before grid disposal (phase $(\mathrm{v})$ in Figure~\ref{fig:MPM_phases}).
Unlike the P2G mapping, the Grid-2-Point (G2P) mapping must always take place, as it is the MPs' role to store information necessary to run the analysis through different steps. 
Under the same assumptions considered in the above paragraph for the P2G mapping (i.e., quasi-static simulations involving an incompressible fluid), the G2P mapping is not required to conserve physical quantities of interest for the same reasons outlined above. 
Hence, the shape functions can be straightforwardly used to map the values of the unknowns from the grid nodes to the MPs, i.e.,
\begin{gather}
{^{mp}\boldsymbol{u}_{n+1}} 
= 
{^{mp}\boldsymbol{u}_{n}}
+
\mathrm{N}^{h/2}_a \left( ^{mp}\boldsymbol{x}_{n+1}  \right) \, \Delta \mathbf{u}_{a};
\\
{^{mp}\boldsymbol{p}_{n+1}} 
=
\mathrm{N}^{h}_a \left( ^{mp}\boldsymbol{x}_{n+1}  \right) \, \left( \mathbf{p}_{n+1} \right)_{a},
\end{gather}
where $\Delta \mathbf{u}_{a}$ is the difference in the time-step of the grid displacements, i.e., $\Delta \left( \bullet \right) \coloneqq \left( \bullet \right)_{n+1} - \left( \bullet \right)_n$. 
The above updates exploit the considerations drawn so far: on the one hand, the initial displacement on the grid $\left( \mathbf{u}_{a} \right)_{n}$ is not reconstructed with the P2G mapping (the unknown of the system of primary equations can directly be $\Delta \mathbf{u}_{a}$); on the other, the pressure field can be re-written both at the MP and grid level because it is not related to an energy measure.

While the update of the mixed volume for the sMPM follows the standard procedure for volume update using the Jacobian, the length of the characteristic function for the GIMPM requires an update which can differ for the different Cartesian directions. This work follows the \emph{corner update procedure} proposed by Coombs \emph{et al.}~\cite{coombs2020lagrangian}. 
The reader is referred to this reference for further details on its implementation, and the exhibited advantages over the other techniques. 

\section{Numerical examples}
\label{sec:numerical_examples}
The model detailed in Section~\ref{sec:mixed_theory}, formulated in a weak sense in~\ref{sec:continuum_formulation}, and discretised in~\ref{sec:MPM_discretisation} is below investigated via three different examples, each exploring a particular feature of the proposed formulation. 
The implementation of the outlined model has been carried out in an extended version of AMPLE~\cite{coombs2020ample}.

\subsection{Terzaghi mono-dimensional consolidation}
\label{subsec:Terzaghi_consolidation}
\paragraph{Example scope}
The numerical investigation provided by the Terzaghi mono-dimensional consolidation (see, for instance,~\cite{terzaghi1943theoretical}) delivers a two-fold goal.
A comparison considering two meshes is made between the Qk$_{SD}$-Qk meshes (specifically using sMPM, i.e., Q1$_{SD}$-Q1) and the Polynomial Pressure Projection (PPP) with a single mesh. 
The PPP is one of the most widely adopted stabilisation techniques in the literature for $\boldsymbol{u}-p^{(f)}$ (see White and Borja~\cite{white2008stabilized} for the FEM, and Zhao and Choo~\cite{zhao2020stabilized} for its adaptation to the MPM).
The results on the two meshes are compared near zero consolidation time versus an established critical time for $\boldsymbol{u}-p^{(f)}$ formulations (see Vermeer~\cite{vermeer1981accuracy}). 
This time sets that the coarser simulation shows oscillating pressure, while the finer analysis does not. While the simulations with Q1$_{SD}$-Q1 meshes agree with these predictions, it is shown that the analyses employing the PPP exhibits pressure instability for both meshes. 
Furthermore, in the case of the coarser unstable simulation, it is investigated how this instability dissipates over time for the simulation using Q1$_{SD}$-Q1 meshes.


The analytical values of pressure are well-known in the small-strain regime and are given by (see, for instance,~\cite{terzaghi1943theoretical}) 
\begin{multline} 
\label{eq:Terzaghi}
P(Z,T) = \sum_{m=0}^{\infty} \frac{2}{M} \, \sin \left( M  Z \right) \exp \left( - M^2 \, T \right), 
\\
\hbox{with} \  M = \frac{\pi}{2} \left( 2 m + 1 \right) \quad \text{and} \ m = 0, 1, \dots \in \mathbb{N}_0,
\end{multline}
where the adimensional quantities in the above equation are defined as follows
\begin{equation}
P \coloneqq \frac{p^{(f)}}{w}; 
\qquad Z \coloneqq \frac{z}{H}; 
\qquad T \coloneqq \frac{c_v}{H^2} \, t,
\end{equation}
with $w$ being the magnitude of the overburden and $H$ the height of the column.
The relationship between the coefficient of consolidation $c_v$ and other hydro-mechanical parameters is given by
\begin{equation}
c_v = \frac{\left(\frac{K}{n_0} + \frac{4}{3} G \right)}{\rho^{(f)} \, g} \, \upkappa,
\end{equation}
where $\frac{K}{n_0}$ is the initial tangent modulus of the improved Hencky material described by Eq.~\eqref{eq:strain_energy}.

\begin{table}
\begin{minipage}[b][8.5cm][c]{0.485\linewidth}
\caption{Summary of the parameters considered in the analysis of the Terzaghi mono-dimensional consolidation.}
\label{table:Terzaghi_setting}
\centering
\begin{threeparttable} 
{\scriptsize{
\begin{tabular}{ccc} 
\headrow
\multicolumn{3}{c}{Parameter Settings}
\\ \hiderowcolors 
\hline
\multirow{3}{*}{\splitcellc{Material \\ Parameters}}
& $\frac{K}{n_0}, \, G$ & $1 \cdot 10^6$ Pa, $ 6 \ \cdot 10^5$ Pa\\
& $\upkappa_0$ & $10^{-7}$ m s$^{-1}$\\
& $ \rho^{(f)}$ & 1000 kg/m$^3$\\
\hline 
\multirow{2}{*}{\splitcellc{Geometry and \\ loading}}
& $H$ & $1$ m \\
& $w$ & $100$ kPa \\
\hline
\end{tabular}
}
}
\end{threeparttable}
\end{minipage}
\hfill
\begin{minipage}[b][7.5cm][c]{0.485\textwidth}
 \centering
  \includegraphics[height=0.3\textheight]{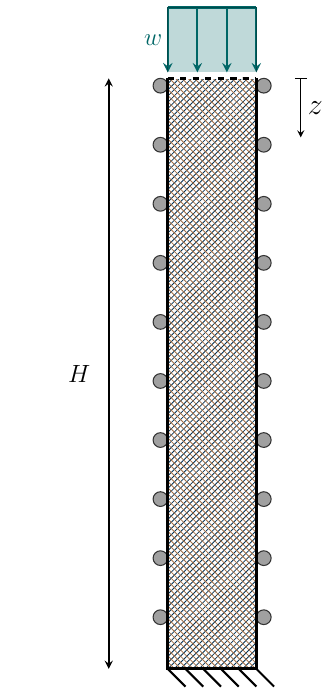}
  \captionof{figure}{Illustration of the Terzaghi mono-dimensional problem. Permeable surfaces are designed by the dashed line.}
  \label{subfig:Terzaghi_column}
\end{minipage} 
\end{table}

\paragraph{Setup}

The hydro-mechanical parameters are the same employed by Zhao and Choo~\cite{zhao2020stabilized} and are reported in Table~\ref{table:Terzaghi_setting}. Owing to the assumed parameters, it is expected that the numerical values, computed in the finite strain context, can reproduce the analytical results set in the small strain theory. 
Two MPs per direction are initially positioned for each element of the finer mesh. 
Owing to this initial discretisation and given that small displacements are expected, the small-cut issue detailed in Section~\ref{subsec:small-cut_issue} does not occur in this example. 
Thus, the ghost stabilisation included in Eqs.~\eqref{eq:eq_equation_discrete_form_QaSD-Qa} and~\eqref{eq:fluid_mass_conservation_discrete_form_QaSD-Qa} are not considered for this problem. 

As mentioned above, two different vertical discretisations have been considered, and the initial selected time-step is $t_0 = 0.1$ s. 
These discretisations have been designed to comply with (in the case of $n_y^{els} = 320$) and violate (in the case of $n_y^{els} = 160$) at the first time-step the time value prescribed by Vermeer and Verruijt~\cite{vermeer1981accuracy} for implicit $\boldsymbol{u}-p^{(f)}$ formulations, this being
\begin{equation}
\label{eq:stable_time-step}
t \geq \frac{\Delta h_y^2}{6 \, c_v}  
\begin{cases}
\approx 0.09 \, \text{s} \qquad \text{if} \ n^{els}_y = 320; 
\\
\approx 0.36 \, \text{s} \qquad \text{if} \ n^{els}_y = 160. 
\end{cases}
\end{equation}
To avoid unnecessary computational time and still provide a good description of the consolidation phenomenon, the time-step partition observes the geometric series to complete the consolidation process (i.e., $T = 1$)
\begin{equation}
\sum_{p=0}^{n-1} c^{p} t_0 \approx \frac{H^2}{c_v}, 
\end{equation}
where the selected common ratio $c = 1.01673$ has been selected to run the simulation in $n = 550$ time-steps. 

The analyses have been run using sMPM shape functions (Q1$_{SD}$-Q1 meshes). 
However, given that small displacements are expected, no appreciable difference could be noticed if the GIMPM shape function should have been employed.

The parameter for the PPP vector is $\frac{1}{2 G}$, which is widely adopted in the literature (see Zhao and Choo~\cite{zhao2020stabilized} and references therein).

\paragraph{Results discussion}
\begin{figure}
\centering
\begin{subfigure}[t]{0.45\textwidth}
\includegraphics[width=\textwidth,trim= 0cm 0.25cm 0 0]{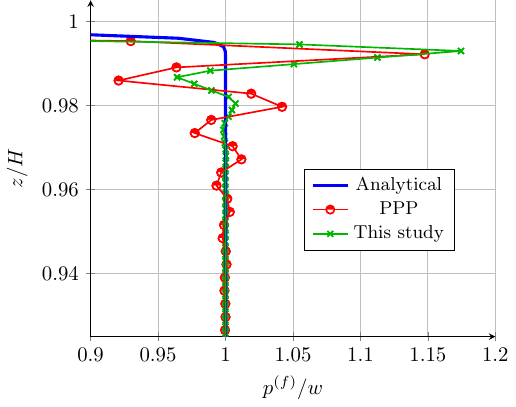}
\caption{Excess pore pressure isochrones comparison at the time $T\approx 1.8 \cdot 10^{-6}$ s for the case with $n^{els}_y = 160$.}
\label{subfig:160els}
\end{subfigure}
\hfill
\begin{subfigure}[t]{0.45\textwidth}
\includegraphics[width=\textwidth,trim= 0cm 0.25cm 0 0]{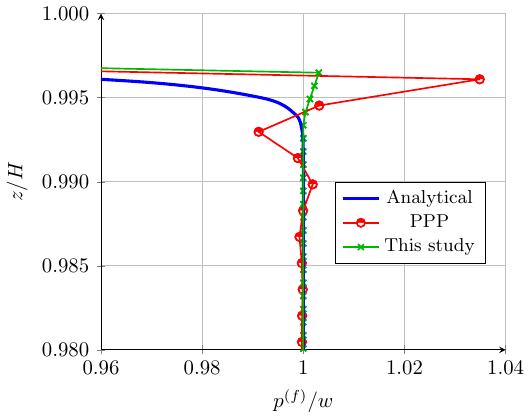}
\caption{Excess pore pressure isochrones comparison at the time $T \approx 1.8 \cdot 10^{-6}$ s for the case with $n^{els}_y = 320$.}
\label{subfig:320els}
\end{subfigure}
\caption{Graphical comparison between analytical formula and numerical solutions obtained with the PPP stabilisation and Q1$_{SD}$-Q1 meshes.}
\label{fig:comparison}
\end{figure}

From a comparison between Figs.~\ref{subfig:160els} and~\ref{subfig:320els}, it is clear that the formulation proposed in this paper (green line) agrees with the prediction of unstable (Fig.~\ref{subfig:160els}) and stable (Fig.\ref{subfig:320els}) time provided by Eq.~\eqref{eq:stable_time-step}. 
Conversely, the PPP (red line) exhibits unstable behaviour even for the predicted stable discretisation with $n_y^{els}=320$. 
Moreover, it can be seen that the instability caused by the PPP propagates through more vertical elements than its Q1$_{SD}$-Q1 counterpart for both situations. 
It must be observed that these results comply with the literature: Preisig and Prevost~\cite{preisig2011stabilization} stated that the PPP is unable to remove pressure oscillation near the draining boundary, and Monforte \emph{et al.}~\cite{monforte2019low}, adopting a $\boldsymbol{u}-\boldsymbol{w}-p^{(f)}$ formulation, showed pressure oscillating behaviour for a high-frequency wave propagation problem. 
These studies confirm that the PPP was not designed to stabilise poro-mechanical problems in which the pressure field belongs to $H^1 (\omega)$: the original application of the PPP (see Dohrmann and Bochev~\cite{dohrmann2004stabilized}) was Stokes problems, which requires the pressure field to be $H^0 (\omega)$. 
The peaks shown in Figs.~\ref{subfig:160els} and~\ref{subfig:320els} (and confirmed by the literature) correspond with sharp pressure gradients (i.e., with the phenomenon progressively passing from undrained to drained), which the PPP was not designed to stabilise.

\begin{figure}[!h]
\centering
\includegraphics[width=0.8\textwidth, trim=0.5cm 0.45cm 0cm 0cm]{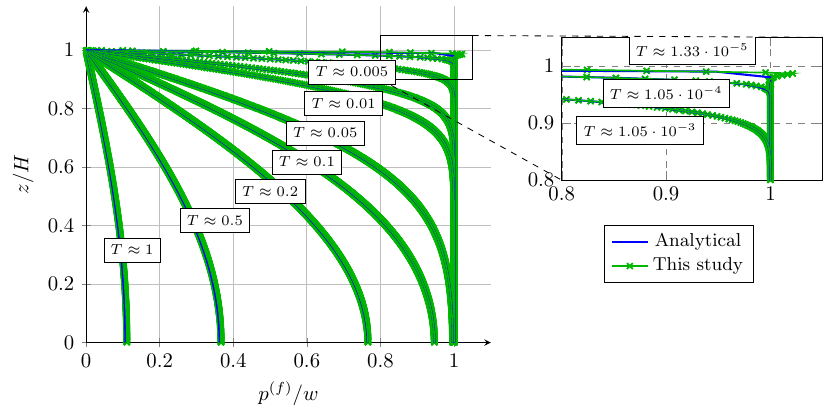}
\caption{Excess pore pressure isochrones for the case with $n^{els}_y = 160$.}
\label{fig:Terzaghi_consolidation}
\end{figure}

Fig.~\ref{fig:Terzaghi_consolidation} (particularly its magnification) investigates the persistence of the spurious peak for the Q1$_{SD}$-Q1 formulation in the unstable case (i.e., $n^{els}_y = 160$). 
The numerical results agree with Eq.~\eqref{eq:stable_time-step}, which prescribes a stable time bigger than $\approx 0.36$ s. 
For the adimensional time $T \approx 1.33 \cdot 10^{-5}$ (corresponding to $\approx 0.74$ s) the peak is practically extinguished. 
Its pressure value is $2.14\%$ bigger than its analytical one, which is significantly lower than the one exhibited in Fig.~\ref{subfig:160els} for $T \approx 1.8 \cdot 10^{-6}$ (corresponding to $\approx 0.1$ s), this being $17.5\%$ bigger compared to the analytical value.
For the $T \approx 1.05 \cdot 10^{-4}$ (corresponding to $\approx 5.83$ s), there is no appreciable difference between the numerical and the analytical solution.
This overlap between the numerical investigation and the analytical solution for the exceed pore pressure isochrones continues for the rest of the simulation, until consolidation has entirely taken place. 

\subsection{Investigation of the face ghost penalty}
\label{subsec:investigation_ghost}

\begin{figure}[!h]
\centering
\includegraphics[width=0.6\textwidth, trim=1.1cm 1cm 0cm 0cm]{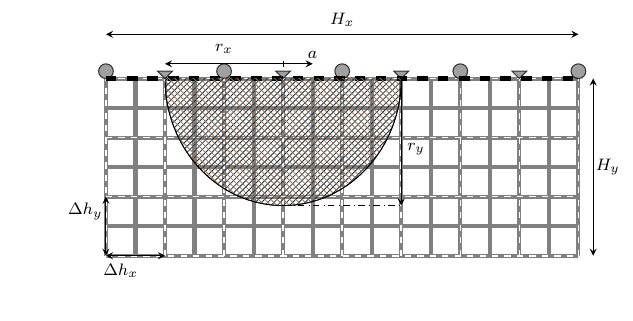}
\caption{Geometry and initial setup for the considered ellipse. 
Each initial setup consists in seeding the MPs discretising the ellipse in a different position, described the coordinate $a$, while keeping the meshes (finer in grey, coarser dashed in white) fixed. 
New MP-setups are generated for different values of $a$.
Drained conditions are applied at the top of the mesh (dashed black line).}
\label{fig:moving_ellipse}
\end{figure}

\paragraph{Example scope}
While the previous example focused on the stabilising effect given by Q1$_{SD}$-Q1 elements, this numerical exercise (illustrated in Figure~\ref{fig:moving_ellipse}) is particularly aimed at investigating the beneficial role of the ghost penalties when applied to the matrices $\mathbf{A}$ and $\mathbf{C}$ appearing in the Jacobian matrix. 
If adopting (as is the case in this example) Qk$_{SD}$-Qk meshes, the drained case requires coercivity of both these bilinear forms associated with the above matrices to guarantee solvability. 

The effectiveness of the face ghost stabilisation applied to Eqs.~\eqref{eq:eq_equation_discrete_form_QaSD-Qa} and~\eqref{eq:fluid_mass_conservation_discrete_form_QaSD-Qa} is investigated for different values of the parameters $\gamma^A$ and $\gamma^C$. 
The spanned values of these parameters have been chosen to be not too small, resulting in insignificant effects, nor too high, leading to potential locking phenomena (see Badia et al.~\cite{badia2022linking}).
This example bears similarities with other numerical tests proposed by Coombs~\cite{coombs2023ghost} in the MPM context and by Kothari and Krause~\cite{kothari2022generalized} and Sticko \emph{et al.}~\cite{sticko2020high} for the unfitted FEM.

\begin{table}[]
\centering
\begin{threeparttable} 
\begin{tabular}{ccc} 
\headrow
\multicolumn{3}{c}{Parameter Settings}
\\ \hiderowcolors 
\hline
\multirow{3}{*}{\splitcellc{Material \\ Parameters}}
& $\frac{K}{n_0}, \, G$ & $1$ Pa, $1$ Pa\\
& $\upkappa_0$ & $1$ m s$^{-1}$\\
& $ \rho^{(f)}$ & $1$ kg/m$^3$\\
\hline 
\multirow{2}{*}{Geometry}
& $H_x, \, H_y$ & $8$ m, $3$ m\\
& $r_x, \, r_y$ & $2$ m, $2.15$ m  \\
\hline
Grid
&
$\Delta h_x, \, \Delta h_y$ & $1$ m, $1$ m \\
\hline
\end{tabular}
\caption{Summary of the parameters considered in the investigation of the ghost penalty example.}
\label{table:moving_ellipse}
\end{threeparttable}
\end{table}

\paragraph{Setup}
The half of an ellipse illustrated in Fig.~\ref{fig:moving_ellipse} is made of an ideal porous material, whose parameters are listed in Table~\ref{table:moving_ellipse}. 
Since the scope of this investigation lies in the examination of the small-cut issue (and its remedy), the geometrical aspect, i.e., the overlap between the physical domain and shape functions' stencil, is examined. 
This justifies two aspects: on the one hand, the assumption of ideal, non-descriptive, hydro-mechanical parameters; on the other, the use of a geometry such an ellipse. 
This form combines a periodic behaviour along the coordinate $a$ (see Fig.~\ref{fig:moving_ellipse}) and aggravates the small-cut problem at the bottom limit of the ellipse. 

Three MPs per direction are initially equally distributed per each element of the finer mesh. 
The MPs lying outside the analytical shape described by the half of the ellipse are then removed. 
This setup is repeated for different values of the coordinate $a$, which is varied progressively for each simulation. 
For each of these discretisations, the submatrices $\mathbf{A}$ and $\mathbf{C}$ in the Jacobian are assembled, so that the effect of the ghost stabilisation can be assessed for different physical MP-based domain-mesh interactions.

\begin{figure}
\centering
\begin{subfigure}[t]{0.45\textwidth}
\includegraphics[width=\textwidth,trim= 0cm 0.5cm 0 0]{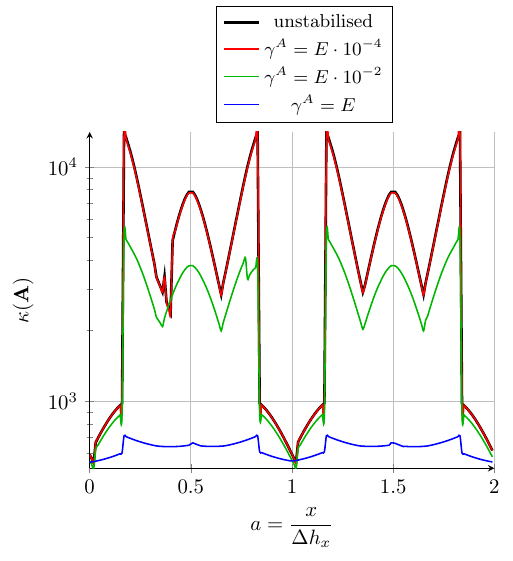}
\caption{Condition numbers for the $\mathbf{A}$ sub-matrix for translating domains, sMPM.}
\label{subfig:MPM_A}
\end{subfigure}
\hfill
\begin{subfigure}[t]{0.45\textwidth}
\includegraphics[width=\textwidth,trim= 0cm 0.5cm 0 0]{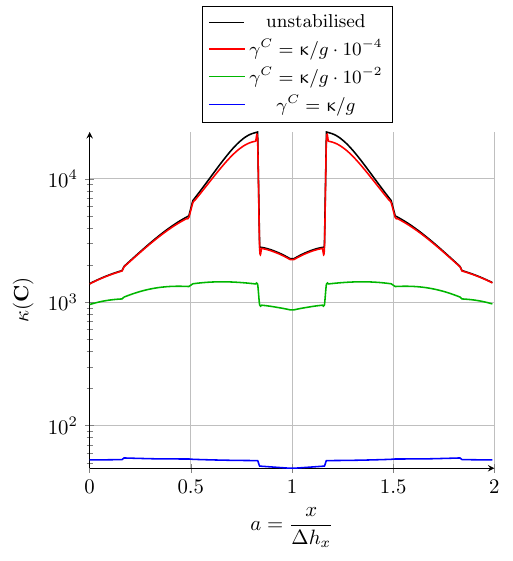}
\caption{Condition numbers for the $\mathbf{C}$ sub-matrix for translating domains, sMPM.}
\label{subfig:MPM_C}
\end{subfigure}
\\
\vspace{0.2cm}
\begin{subfigure}[t]{0.45\textwidth}
\includegraphics[width=\textwidth,trim= 0cm 0.5cm 0 0]{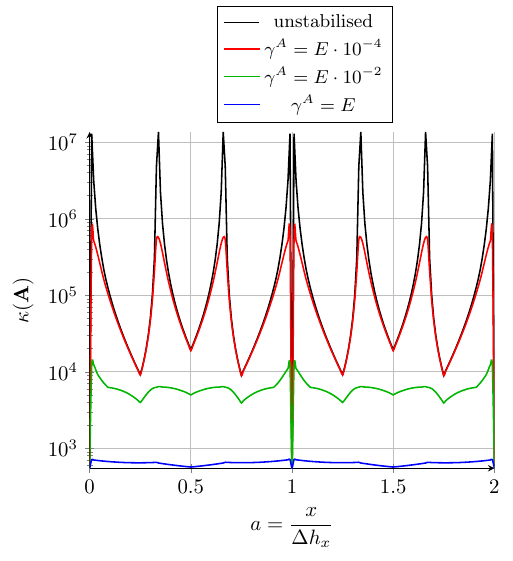}
\caption{Condition numbers for the $\mathbf{A}$ sub-matrix for translating domains, GIMPM.}
\label{subfig:GIMPM_A}
\end{subfigure}
\hfill
\begin{subfigure}[t]{0.45\textwidth}
\includegraphics[width=\textwidth,trim= 0cm 0.5cm 0 0]{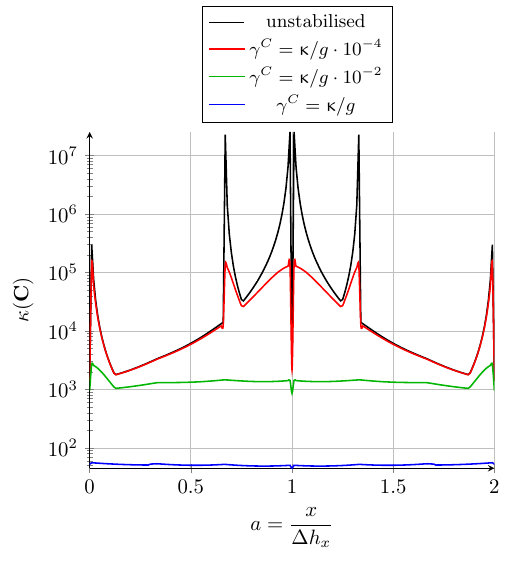}
\caption{Condition numbers for the $\mathbf{C}$ sub-matrix for translating domains, GIMPM.}
\label{subfig:GIMPM_C}
\end{subfigure}
\caption{Effects of the applied ghost penalty parameters on the condition number of submatrices $\mathbf{A}$ (left column) and $\mathbf{C}$ (right column) in the case of sMPM (top row) and GIMPM. (bottom row).}
\label{fig:condition_numbers}
\end{figure}

It must be noted that, for the sMPM, it is not possible to assign an element to a MP when this lie precisely on an element boundary. 
Since this occurs for the finer mesh when $a = \frac{\Delta h_{x}}{4} * p$, with $p = 1, 2, \dots \in \mathbb{N}$, these sampling location are avoided for the sMPM.
This choice explains the rupture of periodic behaviour which is particularly evident in Figure~\ref{subfig:MPM_A}. 
Nonetheless, it must be underlined that this situation is due to the specificity of the setup and is highly unlikely to occur in standard simulations.

\paragraph{Results discussion}

As it can be appreciated from Fig.~\ref{fig:condition_numbers}, there is considerable difference in the condition numbers (denoted by $\upkappa(\bullet)$, with $\bullet$ being the considered submatrix) for unstabilised sMPM and GIMPM, with the latter showing on average three orders of magnitude higher condition number for the same setup. This behaviour is due to the shape functions stencil, which is more extended for the GIMPM.

Fig.~\ref{subfig:MPM_A} highlights how the lowest value of ghost penalty ($\gamma^A = E \cdot 10^{-4}$) does not contribute significantly to the submatrix $\mathbf{A}$ condition number in the case of sMPM. 
Diversely, in the case of the GIMPM in Fig.~\ref{subfig:GIMPM_A}, the effect of the ghost penalty for $\gamma^A = E \cdot 10^{-4}$ reduces the peak values by approximately an order of magnitude, even though the general pattern follows similarly that of the unstabilised version. 
These peaks are more levelled for the case of $\gamma^A = E \cdot 10^{-2}$ and entirely smoothed for the case of $\gamma^A = E$, which exhibits condition numbers approximately 4.5 order of magnitude lower than the peaks of the unstabilised version. 
Back to Fig.~\ref{subfig:MPM_A} for the sMPM, it can be seen how $\gamma^A = E$ leads to overall smooth behaviour, with the condition number being approximately 1.5 orders of magnitude below the unstabilised peaks. 
The value $\gamma^A = E \cdot 10^{-2}$ contributes less significantly to stabilising the sMPM than its GIMPM counterpart.
Similar trends can be seen for the condition number of the submatrix $\mathbf{C}$ in the case of sMPM (Fig.~\ref{subfig:MPM_C}) and GIMPM (Fig.~\ref{subfig:GIMPM_C}). 
This behaviour is interesting, as the submatrix $\mathbf{C}$ is entirely assembled using the coarser grid of length $h$. 
Hence, it highlights how coarser discretisation, which should be less prone to the small-cut issue when using the same numbers of MPs, is not free from this problem and strengthens the need for a stabilisation such as the one included in this work.

Owing to the previously mentioned lack of rigorous analysis for the MPM, the condition numbers of the submatrices $\mathbf{A}$ and $\mathbf{C}$, even if they benefit from the ghost stabilisation, are still dependent to some degree on the intersection between the shape functions' stencil and the MP-based physical domain.
This aspect diversifies the MPM from the unfitted FEM, where, once stabilised, the domain's cuts do not play a role in the condition number (as proved in Burman~\cite{burman2010ghost} and numerically demonstrated in Kothari and Krause~\cite{kothari2022generalized}). This difference amplifies the importance of carefully choosing the correct value of the stabilisation parameter, as some selections, even though effective, might not give the desired limitation on the condition number.

\subsection{Flexible strip foundation}
\label{subsec:flexible_strip_footing}

\begin{table}
\centering
\centering
\begin{threeparttable} 
{\footnotesize{
\begin{tabular}{ccc} 
\headrow
\multicolumn{3}{c}{Parameter Settings}
\\ \hiderowcolors 
\hline
\multirow{5}{*}{\splitcellc{Material \\ Parameters}}
& $K, \, G$ & $1.062 \cdot 10^8$ Pa, $ 9.8 \cdot 10^7$ Pa\\
& $\alpha, \, \gamma, \, p_c$ & $0.6, \, 0.9, \, 7.5 \cdot 10^7$ Pa\\
& $\upkappa_0$ & $2.07 \cdot 10^{-3}$ m s$^{-1}$\\
& $n_0$ & 0.5\\
& $ \rho^{(f)}$ & $1000$ kg/m$^3$\\
\hline 
\multirow{2}{*}{\splitcellc{Geometry and \\ loading}}
& $L, \, H, \, l$ & $20$ m, $16$ m, $8$ m\\
& $w$ & $2.5 \cdot 10^8$ Pa \\
\hline
\multirow{2}{*}{\splitcellc{Time and \\ time-step}}
&
T & $10$ s \\
&
$\Delta t$ & $0.5$ s \\
\hline
\multirow{3}{*}{\splitcellc{Grid \\ and MPs}}
&
$\Delta h_x, \, \Delta h_y$ & $0.25$ m, $0.25$ m
\\
&
\multirow{2}{*}{\splitcellc{MPs per element \\ (finer mesh)}} 
& \multirow{2}{*}{4} 
\\
\\
\hline
\end{tabular}
\caption{Summary of the parameters considered in the analysis of the flexible strip foundation.}
\label{table:flexible_footing}
}
}
\end{threeparttable}
\end{table}

\begin{figure}
\centering
\includegraphics[width=0.45\textwidth,trim= 0.5cm 0 0 0]{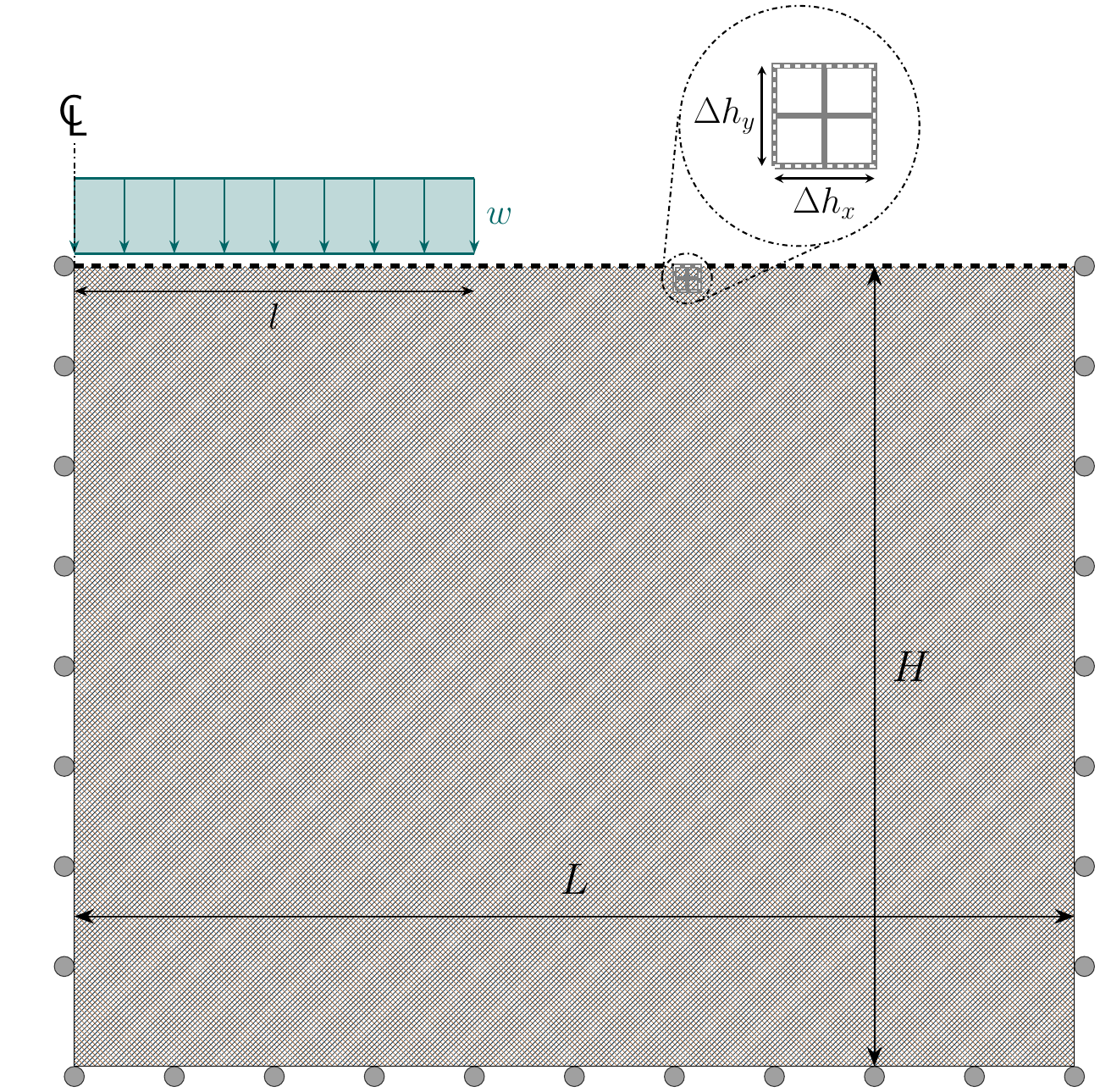}
\caption{Illustration of the flexible strip foundation problem. Permeable surfaces are designed by the dashed line.}
\label{subfig:foundation_geometry}
\end{figure}

\paragraph{Example scope}

The scope of this numerical test is to support the rationale provided in Section~\ref{subsubsec:GIMPM} regarding the stability of the different combinations of polynomial functions employed by GIMPM formulation with Qk$_{SD}$-Qk meshes. 
To demonstrate that the pressure field is free from oscillations, a flexible footing loading has been applied at the top of a fully saturated porous material in plane-strain conditions.
Elastic and the elasto-plastic cases have been considered.
Since the application of the loading causes the considered body to displace substantially, it is expected that all the combinations of polynomials that the GIMPM basis functions exhibit will be met through the analysis, even if not uniformly or simultaneously.

\paragraph{Setup}
Hydro-mechanical parameters are reported in Table~\ref{table:flexible_footing}, and the considered elastic material shares similarities with the one presented in Armero~\cite{armero1999formulation} in a very similar example. 
The modified Cam clay $\alpha-\gamma$ yield surface proposed by Collins and Hilder~\cite{collins2002theoretical} has been considered for the elasto-plastic case (see, for instance,~\cite{coombs2011algorithmic} for implementation details).

To apply the non-conforming pressure BCs, a penalty value of $\gamma^{pen} = 5 \cdot 10^6 / \upkappa$ was employed, where the hydraulic conductivity appears at the denominator for a matter of consistency with the physical units.
The load has been linearly increased from zero to the value shown in Table~\ref{table:flexible_footing}. The elastic case can complete the simulation in 20 time-step, while the plastic case reaches the limit load in correspondence of the $4^{th}$ time-step, where it is stopped.

\begin{figure}
\centering
\begin{subfigure}[t]{0.45\textwidth}
\includegraphics[width=\textwidth,trim= 0cm 0cm 0 0]{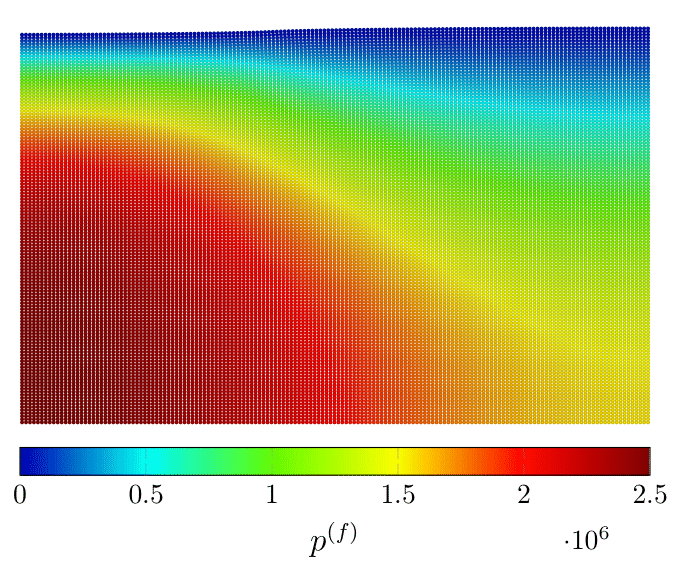}
\caption{Contours of the fluid pressure field, $1^{st}$ time-step.}
\label{subfig:1}
\end{subfigure}
\\
\begin{subfigure}[t]{0.45\textwidth}
\includegraphics[width=\textwidth,trim= 0cm 0cm 0 0]{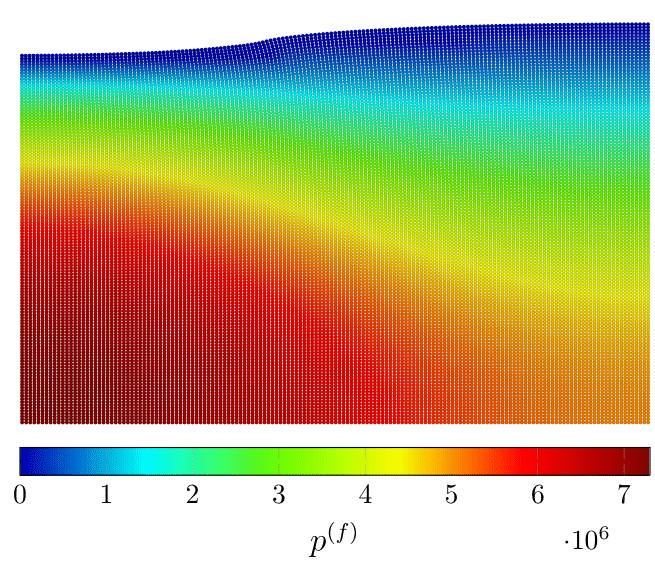}
\caption{Contours of the fluid pressure field, $5^{th}$ time-step.}
\label{subfig:5}
\end{subfigure}
\hfill
\begin{subfigure}[t]{0.45\textwidth}
\includegraphics[width=\textwidth,trim= 0cm 0cm 0 0]{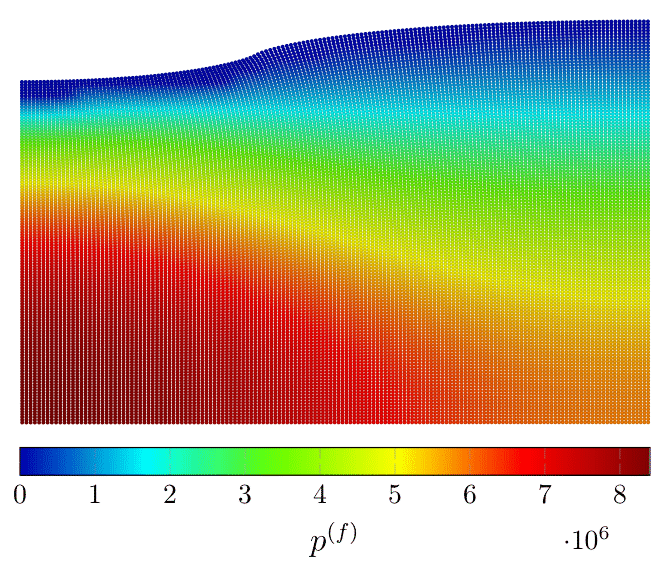}
\caption{Contours of the fluid pressure field, $10^{th}$ time-step.}
\label{subfig:10}
\end{subfigure}
\\
\begin{subfigure}[t]{0.45\textwidth}
\includegraphics[width=\textwidth,trim= 0cm 0cm 0 0]{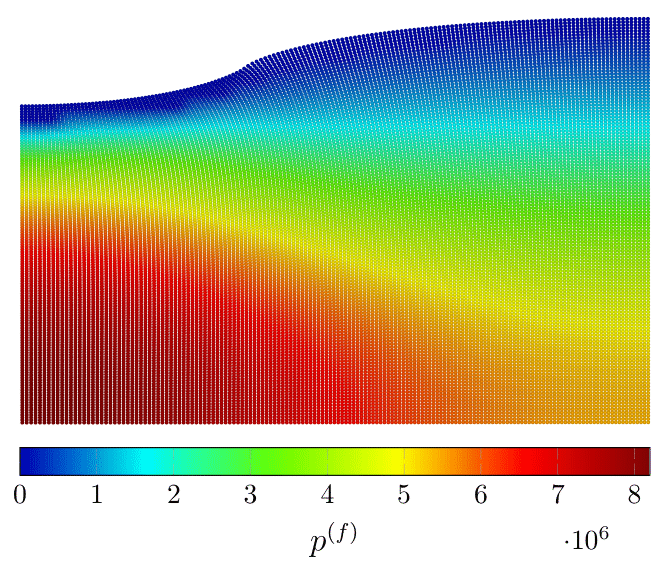}
\caption{Contours of the fluid pressure field, $15^{th}$ time-step.}
\label{subfig:15}
\end{subfigure}
\hfill
\begin{subfigure}[t]{0.45\textwidth}
\includegraphics[width=\textwidth,trim= 0cm 0cm 0 0]{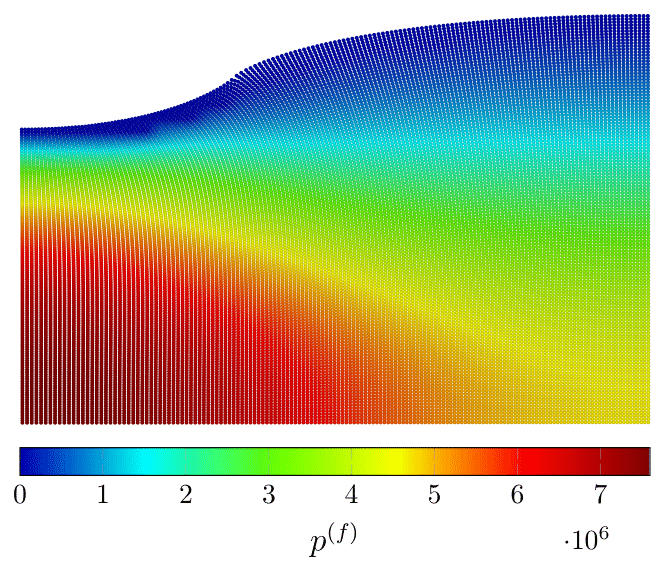}
\caption{Contours of the fluid pressure field, $20^{th}$ time-step.}
\label{subfig:20}
\end{subfigure}
\caption{Contours of the pressure field across the simulation for the flexible footing loading, elastic case.}
\label{fig:contours}
\end{figure}

Given the low number of MPs per element and the substantial displacements expected in this simulation, the ghost penalisation have been included. 
Based on the observations made in Example~\ref{subsec:investigation_ghost}, the selected parameter for the penalty value are $\gamma^A = E \cdot 10^{-1}$ and $\gamma^C = \frac{\upkappa}{g} \cdot 10^{-1}$.

\paragraph{Results discussion}

\begin{figure}
\centering
\begin{subfigure}[t]{0.45\textwidth}
\includegraphics[width=\textwidth,trim= 0cm 0cm 0 0]{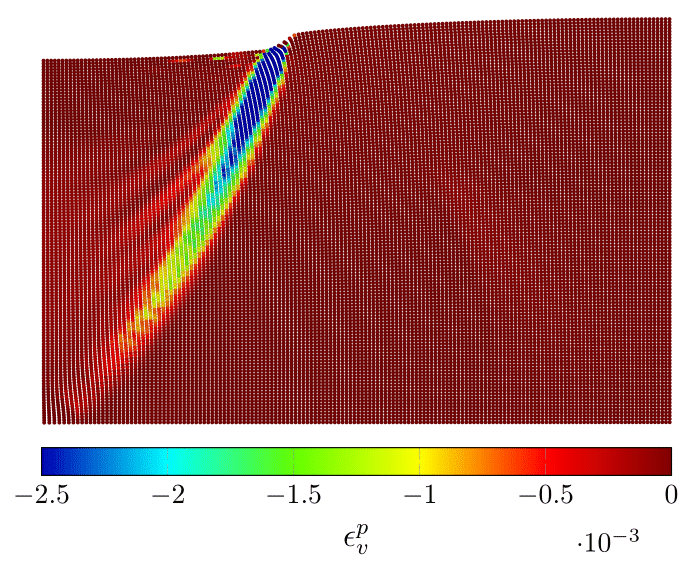}
\caption{Contours of the plastic volumetric part of the logarithmic strain, last time-step.}
\label{subfig:plastic_strain}
\end{subfigure}
\hfill
\begin{subfigure}[t]{0.45\textwidth}
\includegraphics[width=\textwidth,trim= 0cm 0cm 0 0]{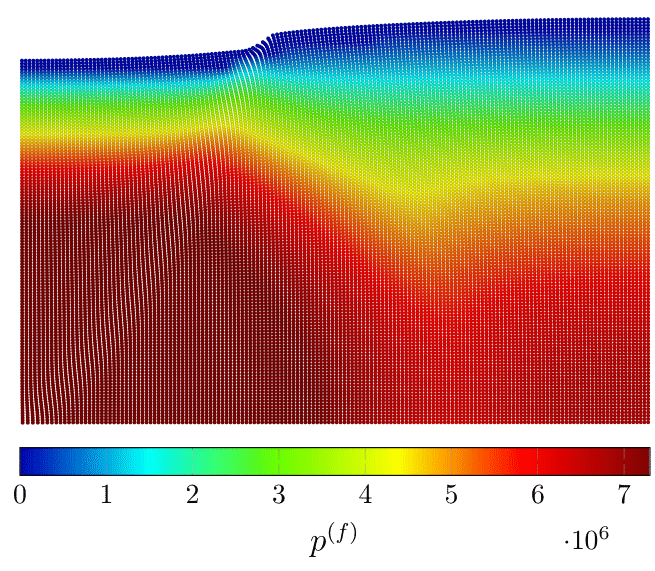}
\caption{Contours of the fluid pressure field, last time-step.}
\label{subfig:pressure}
\end{subfigure}
\caption{Contours of different fields at the last time-step for the flexible footing loading, plastic case.}
\label{fig:contours_elastoplasticity}
\end{figure}

As can be seen from the contours shown throughout the analysis in Fig.~\ref{fig:contours} for the elastic case, the fluid pressure does not present oscillations typical of violating the inf-sup condition. 
Even though showing spatially different distributions of pressure, these values are equally smooth for the final step of the elasto-plastic case shown in Fig.~\ref{subfig:pressure}.
As Fig.~\ref{subfig:plastic_strain} shows, plastic strains arise in correspondence with the footing angle, forming the expected wedge shape. 
This shape is also visible for the pressure field in Fig.~\ref{subfig:pressure}, where the GIMPM Qk$_{SD}$-Qk meshes handle this horizontal gradient without interruptions or oscillations.
Due to the continuity demonstrated for the different fields and situations in Figs.~\ref{fig:contours} and~\ref{fig:contours_elastoplasticity}, the chosen ghost parameters have proved suitable for stabilising the submatrices' conditioning number and achieving the continuous expected results. 

Overall, the exhibited smooth values of pressure substantiate the explanation provided in Section~\ref{subsubsec:GIMPM} for the GIMPM shape functions, showing that these described Qk$_{SD}$-Qk stable meshes deliver oscillating-free solutions.

\section{Conclusions}
This manuscript has highlighted two sources of instability that can arise for the MPM in the case of mixed formulations and proposed a new approach that overcomes both issues for coupled (solid-fluid) problems.
The intrinsic nature of the MPM as an unfitted method and its resulting small-cut issue can lead to ill-conditioned matrices. 
For mixed formulations, this small-cut issue affects both the submatrices appearing on the main diagonal of the Jacobian (i.e., $\mathbf{A}$ and $\mathbf{C}$), especially in drained conditions. 
The treatment of this instability has been the use of the face ghost penalty method on both the displacement and fluid pressure meshes, which, even though it does not guarantee coercivity for the MPM, limits the condition number of these submatrices, thus permitting the inversion of the Jacobian matrix.

Furthermore, in nearly undrained conditions, the choices of the spaces of test and trial functions relative to displacement and pressure can violate the inf-sup condition. 
This work has adapted the use of Qk$_{SD}$-Qk elements to the MPM. The resulting overlapping meshes, while maintaining low-order of shape functions (as is the case with the sMPM and the GIMPM), are stable by design. This is in contrast to the widely adopted PPP formulation, which is unable to stabilise poro-mechanical problems, especially near drained boundaries, due to the $H^1 (\omega)$ nature of the pressure field \cite{preisig2011stabilization}. 

The key contribution offered by this paper is a stable, implicit MPM formulation for large deformation coupled problems including elasto-plastic material behaviour. 
While a rigorous analysis or even the patch test are not available for the MPM, the rationales provided throughout the paper are substantiated via the similarities between the MPM and the unfitted FEM and tested with the considered numerical examples.


\section*{Acknowledgements}
The first author was supported by the Faculty of Science at Durham University. 
The second, fourth and fifth were supported by Engineering and Physical Sciences Research Council, Grant Number EP/W000970/1. The third author was supported by Physical Sciences Research Council, Grant Number EP/T518001/1.
All data created during this research are openly available at \href{https://collections.durham.ac.uk/}{collections.durham.ac.uk/} (specific DOI to be confirmed if/when the paper is accepted).
The research
presented in this article has also benefited from discussions with, and feedback from Ted O’Hare, Bradley Sims, Sam Sutcliffe, and Dr. Mao Ouyang.
For the purpose of open access, the author has applied a Creative Commons Attribution (CC BY) licence to any Author Accepted Manuscript version arising.



\bibliographystyle{elsarticle-num}
\bibliography{sample}


\end{document}